\documentclass[10pt,a4paper]{article}
\usepackage[english]{babel}
\usepackage{epsfig}
\usepackage{amsmath,amsthm}
\usepackage{amssymb}
\usepackage{amsbsy}
\usepackage{amsfonts}
\usepackage{array}
\usepackage{verbatim}
\usepackage{graphicx}
\usepackage{subfigure}
\usepackage{tabularx}
\usepackage{bm}
\usepackage{booktabs}
\usepackage[font=small]{caption}
\usepackage{stmaryrd}
\usepackage{mathtools}
\usepackage{color}
\usepackage{titlesec}
\usepackage[left=2.5cm,right=2.5cm,bottom=3.5cm,top=3.5cm]{geometry}

\titleformat*{\section}{\large\bfseries}
\titleformat*{\subsection}{\bfseries}
\titleformat*{\subsubsection}{\bfseries}

\DeclareMathOperator\erf{erf}

\newcommand{\q}{\mathbf{q}}

\begin{document}

\begin{center}
\textbf{ \Large{A hyperbolic reformulation of the Serre-Green-Naghdi model \\[6pt] for general bottom topographies} }

\vspace{0.5cm}
{C. Bassi$^{(a,b)}$,
 L. Bonaventura $^{(a)}$, S. Busto$^{(b,c)}$, M. Dumbser$^{(c)}$\footnote{Corresponding author.\\ \hspace*{0.25cm} Email addresses: caterina.bassi@polimi.it (C. Bassi), luca.bonaventura@polimi.it (L. Bonaventura), saray.busto@unitn.it (S. Busto), michael.dumbser@unitn.it (M. Dumbser)}}

\vspace{0.2cm}
 {\small
\textit{$^{(a)}$ MOX--Modelling and Scientific Computing, Dipartimento di Matematica, Politecnico di Milano, Piazza Leonardo da Vinci 32, 20133 Milano, Italy}

\textit{$^{(b)}$ Istituto Nazionale di Alta Matematica ``Francesco Severi'', Piazzale Aldo Moro 5, 00185 Roma, Italy}

\textit{$^{(c)}$ Laboratory of Applied Mathematics, DICAM, University of Trento, via Mesiano 77, 38123 Trento, Italy}
}
\end{center}

\vspace{0.4cm}
	
\hrule
\vspace{0.4cm}

\noindent \textbf{Abstract}

\vspace{0.1cm}
\noindent We present a novel hyperbolic reformulation of the Serre-Green-Naghdi (SGN) model for the description of dispersive water waves. Contrarily to the classical Boussinesq-type models, it contains only first order derivatives, thus allowing to overcome the numerical difficulties and the severe time step restrictions arising from higher order terms. The proposed model reduces to the original SGN model when an artificial sound speed tends to infinity. Moreover, it is endowed with an energy conservation law from which the energy conservation law associated with the original SGN model is retrieved when the artificial sound speed goes to infinity. 
The governing partial differential equations are then solved at the aid of high order ADER discontinuous Galerkin finite element schemes. 
The new model has been successfully validated against numerical and experimental results, for both flat and non-flat bottom. For bottom topographies with large variations, the new model proposed in this paper provides more accurate results with respect to the hyperbolic reformulation of the SGN model with the mild bottom approximation recently proposed in \cite{escalante:2018}.


\vspace{0.2cm}

\vspace{0.4cm}

\hrule

\section{Introduction}
In a wide variety of situations, the propagation of water waves can be successfully described  employing the classical shallow water (SW) equations \cite{saintvenant:1871}. However, a serious drawback of the SW equations is their lack of ability to represent dispersive waves, non-hydrostatic effects and solitary wave propagation. 

The evolution of water waves is actually controlled by both nonlinear effects, which cause wave steepening, and dispersion effects, which are responsible for stabilisation. The propagation of solitary waves is an example of the perfect counterbalancing between these two opposite phenomena. These effects are usually associated with two different parameters \cite{kirby:2016,fernandeznieto:2018}. The first one, associated with dispersion, is denoted by $\mu$ and is proportional to the ratio between the characteristic water depth and the characteristic wavelength of the water waves. Meanwhile the second parameter, $\delta$, is a nonlinearity parameter given by the ratio between the wave amplitude and the water depth. In the literature, a wide variety of dispersive models is obtained by building asymptotic expansions with respect to the dispersion parameter $\mu$ and to the nonlinearity parameter $\delta$ and retaining the terms up to a certain order. 
The pioneering work \cite{boussinesq:1872} provides the first example of Boussinesq-type model for flat bottom topography in one space dimension. The model is based on the assumption of weak dispersion and weak nonlinearity. 
A two-dimensional extension to non-flat bottom geometries is provided instead in \cite{peregrine:1967}, where, however, the limitation concerning nonlinear terms is still maintained. 
A fully nonlinear approach, keeping the weak dispersion hypothesis, is presented in \cite{serre:1953} for the one-dimensional case with flat bottom. A two-dimensional extension for arbitrary bottom, the Serre-Green-Naghdi (SGN) model, is proposed for the first time in \cite{green:1976} and successively in \cite{seabrasantos:1987}, where the model is also successfully tested against experimental data. In \cite{cienfuegos:2006} a derivation of the model in \cite{seabrasantos:1987} is provided using asymptotic expansions and assuming irrotational flow. 

A peculiar feature of Boussinesq-type models is that they contain derivatives of order higher than one. As highlighted in \cite{fernandeznieto:2018} for the model proposed in \cite{seabrasantos:1987,cienfuegos:2006}, it is often possible to rewrite higher-order derivatives by employing auxiliary variables obtaining, as a consequence, augmented systems which contain only first order derivatives. First order systems for the description of dispersive waves can also be directly obtained following a depth averaging procedure similar to the one employed for the derivation of the classical SW equations and retaining non-hydrostatic contributions in the vertical momentum equation up to the desired order of accuracy with respect to a parameter on the ratio of water depth and a typical horizontal length scale. Such a procedure is followed e.g. in \cite{saintmarie:2011}.

Up to now, we have presented literature for which the weak dispersion hypothesis is assumed to be valid. However, over the years a great effort has been devoted to the development of models with better dispersion characteristics, in order to improve the range of applicability of the models. A few examples can be found in \cite{madsen:1991,madsen:1992,nwogu:1993,madsen:2003} and in the review \cite{madsen:2010}.
Extensions of dispersive models have also been proposed to represent additional physical phenomena lying  outside the classical formulation for inviscid flows. Among them we recall wave breaking, turbulence, vorticity and wind effects. For a review on this topic see \cite{kirby:2016}. 

From the numerical point of view, the main drawback of Boussinesq-type models is that, contrarily to the classical hyperbolic SW equations, they contain higher order space and mixed space-time derivatives, which are numerically very challenging to deal with, see \cite{escalante:2018,dumbser:2016}, and that introduce severe time step restrictions when explicit time integration methods are employed.  
A possible solution to this problem could be the introduction, when possible, of augmented first order systems, as the one proposed in \cite{fernandeznieto:2018}. Also in this case, however, the hyperbolicity of the SW equations is lost and the solution of an additional \textit{elliptic} equation is required at each time step. 
An alternative approach has been very recently presented in \cite{escalante:2018}, where a hyperbolic reformulation is proposed for the model introduced in \cite{saintmarie:2015} and for the SGN model with the mild bottom assumption (higher order derivatives for the bottom are considered to be negligible). Notice that the idea of introducing a hyperbolic approximation of a non-hyperbolic system comes from the seminal paper \cite{cattaneo:1958}, where second order derivatives in the heat equation are replaced by relaxation terms. 
A similar approach that allows to rewrite the compressible Navier-Stokes equations as an extended hyperbolic relaxation system was recently forwarded in \cite{PeshRom2014,GPRmodel}.  
 Besides, the first hyperbolic reformulation of a dispersive system has been recently derived from variational principles in \cite{favrie:2017} for the SGN system in the flat bottom case. Other related work on hyperbolic models for dispersive water waves can be found, for example, in~\cite{ricchiutoHyp,brocchini}. 

In the present work we adopt a similar approach as the one employed in \cite{escalante:2018}. In particular, we build a new first order hyperbolic reformulation for the SGN system without the mild bottom approximation (in the first order form presented in \cite{fernandeznieto:2018}) by introducing two evolutionary equations; one for the depth averaged non-hydrostatic pressure and another for the non-hydrostatic pressure evaluated at the bottom boundary. The present approach follows the ideas of the method of artificial compressibility employed in the incompressible Navier-Stokes equations context and it is also similar to the hyperbolic divergence cleaning procedure introduced in \cite{munz:2000,dedner:2002} for the Maxwell and magnetohydrodynamics equations. 
Apart from hyperbolicity, the proposed system satisfies an important additional property.   Indeed, an extra energy conservation law holds for the hyperbolic system proposed in this paper, and it reduces to the energy conservation law associated to the original SGN system when the artificial sound speed tends to infinity. Such an extra conservation law is very important in the context of symmetric hyperbolic and thermodynamically compatible systems, see \cite{God1961,FriedLax1971,Rom1998}.  

Regarding the numerical discretization of the derived hyperbolic reformulation, we will employ an explicit Discontinuous Galerkin (DG) method \cite{chavent:1989,cockburn:1989,cockburn:1989b,cockburn:1990,cockburn:1991}. 
The DG method has been applied for the first time to equations containing higher order derivatives in \cite{yan:2002a}, where the LDG method \cite{bassi:1997,cockburn:1998} is used for the resolution of linear dispersive Korteveg-de-Vries (KdV) equations, containing up to third order spatial derivatives. Extensions to linear equations with derivatives up to fifth order and nonlinear dispersive equations are instead presented in \cite{yan:2002b,levy:2004}. Applications of the DG method to the solution of nonlinear Boussinesq-type dispersive equations have been introduced in \cite{eskilsson:2006a,eskilsson:2006b,engsig:2008}. Notice that, as pointed out in \cite{yan:2002a}, the DG method, if associated with explicit time integration schemes, can be effectively applied to equations containing higher order derivatives only in the context of convection dominated problems, due to the severe time step restrictions introduced by higher order terms in dispersion dominated problems. This difficulty is overcome in \cite{dumbser:2016}, where a fully implicit space-time DG method is applied to both linear third order KdV equations and to nonlinear Boussinesq-type equations. For novel residual  distribution (RD) schemes applied to Boussinesq-type equations the reader is referred to \cite{ricchiuto}, while high order accurate WENO schemes for nonlinear, nonhydrostatic water waves were discussed in \cite{Kontos}.  

Due to the \textit{hyperbolic} character of the new system we are going to propose in this paper, the usual CFL condition holds, hence $\Delta t$ is proportional to $\Delta x$, thus avoiding higher powers of $\Delta x$ in the stability condition that are typical from Boussinesq-type equations. In this context, it is clear that standard DG methods are very well suited for the new model.

Attaining high order of accuracy in space is straightforward in the DG framework, while suitable high order time discretizations are still a very active field of research. A successful approach is the use of the already mentioned space-time DG methods,  \cite{Sander2012,Rhebergen2013,spacetimedg1,spacetimedg2,3DSIINS,3DSICNS,KlaijVanDerVegt,TD18LE,BTBD20}. An alternative, that will be followed in this work, are the ADER-DG schemes first put forward in \cite{dumbser_jsc} and generalized to the unified $P_{N}P_{M}$ framework for arbitrary high order accurate finite volume (FV) and DG schemes in \cite{Dumbser2008}. This methodology can also be seen as an extension of classical ADER methods,  \cite{mill,TMN01,toro4,titarevtoro,dumbser_jsc,Gassner2011a}, which avoids the cumbersome Cauchy-Kovalevskaya procedure, resulting in more general algorithms.  
ADER-DG methods have already been successfully applied also to non-conservative hyperbolic systems and geophysical flows in \cite{ADERNC}, which makes them a suitable choice to discretize also the new hyperbolic reformulation of the SGN model without mild bottom approximation proposed in this paper.

The paper is organized as follows. Section \ref{sec:governingequations} provides the description of the new hyperbolic reformulation of the SGN model without the mild bottom approximation. 
Section \ref{sec:numscheme} is devoted to a brief description of the ADER-DG scheme. Numerical results are presented in section \ref{sec:numtest}, while conclusions and perspectives for future work are provided in section \ref{sec:conclusion}. 

\section{Governing equations} \label{sec:governingequations}
As suggested in \cite{fernandeznieto:2018}, the SGN system without the mild bottom approximation can be rewritten as a first order system as

\begin{subequations}\label{eq:gnb2}
\begin{align}
& \partial_t h+ \partial_x(h\overline u) = 0, \label{eq:gnb_mass2} \\
& \partial_t (h\overline u) + \partial_x \left( h{\overline u}^2 +h\overline p\right) +gh\partial_x h + (gh + p_b)\partial_x z_b =0, \label{eq:gnb_momx2}\\
& \partial_t (h\overline w) + \partial_x(h \overline u \ \overline w) = p_b, \label{eq:gnb_momz2}\\
& \partial_t(h\sigma) + \partial_x(h\overline u\sigma)= -6 p_b +12 \overline p  \label{eq:gnb_52}, \\
& \sigma = -h \partial_x \overline u \label{eq:gnb_53a}\\
& \overline w + \frac{1}{2}h\partial_x \overline u - \overline u\partial_x z_b = 0,\label{eq:gnb_42} 
\end{align} 
\end{subequations}
where $h(x,t)$ is the water depth, $\overline u (x,t)$ is the depth averaged horizontal velocity, $p_b(x,t)$ is the non-hydrostatic pressure at the bottom boundary, $\overline p (x,t)$ is the depth averaged non-hydrostatic pressure, $z_b(x)$ is the vertical coordinate of the bottom boundary, $\overline w (x,t)$ is the depth averaged vertical velocity and $\sigma$ is an auxiliary variable equal to $ -h \partial_x \overline u$. Since $z_b$ depends only on $x$, we assume that the bottom boundary can vary in space but is fixed in time.  

\begin{figure}
\centering
\includegraphics[width=0.7\textwidth]{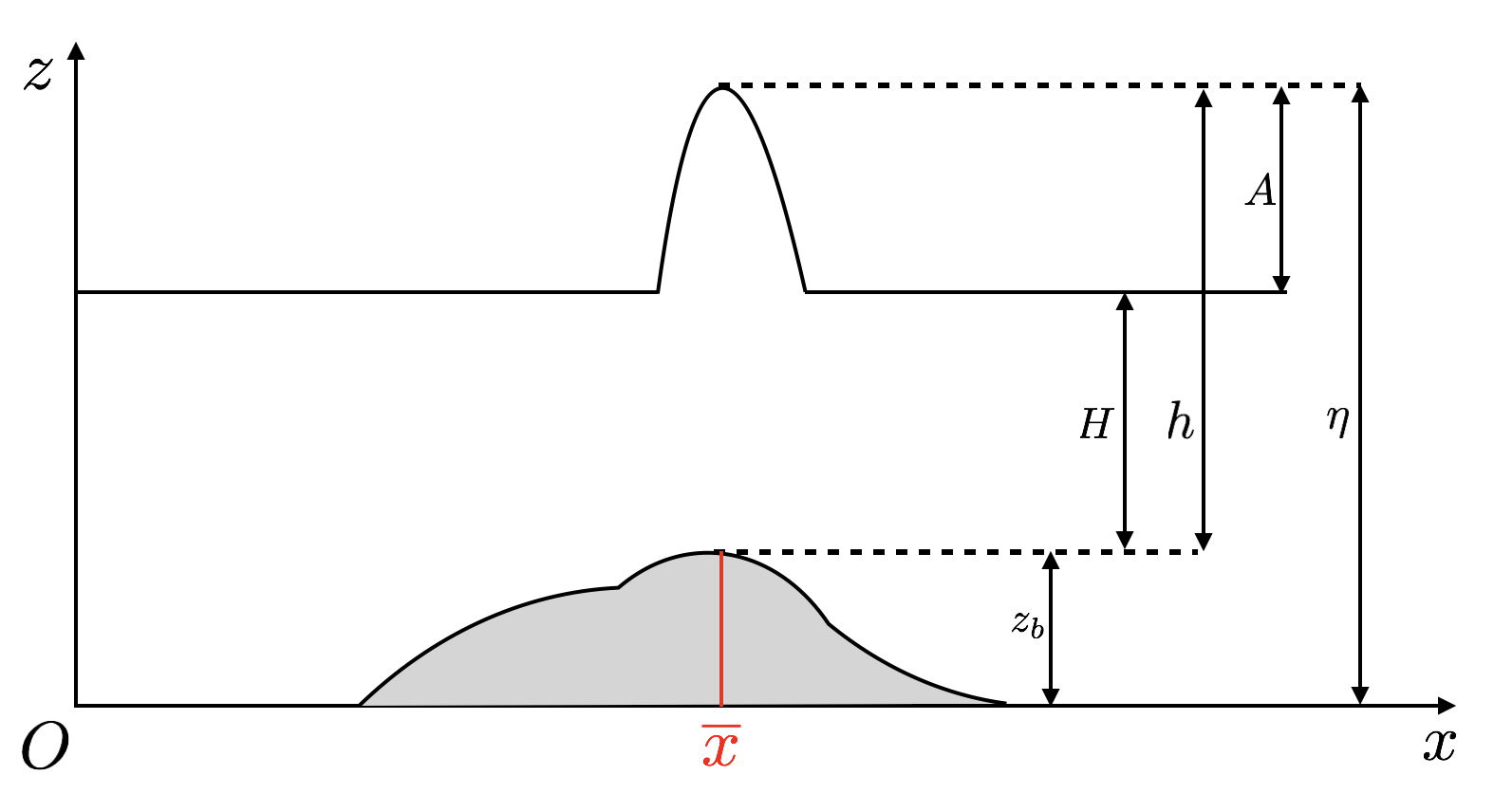}
\caption{Sketch of the shallow water domain. $z_b(x)$ is the vertical coordinate of the bottom boundary, $H(x)$ is the still water depth, $h(x,t)$ is the water depth, $A(x,t)$ is the wave amplitude and $\eta(x,t)$ is the free surface elevation.}
\label{fig:sw_scheme}
\end{figure}

For the sake of clarity, in Figure \ref{fig:sw_scheme}, we represent the coordinate system we employ throughout the paper: besides from $h(x,t)$ and $z_b(x)$ already defined, $\eta(x,t)$, $H(x)$ and $A(x,t)$ are the free surface elevation, the still water depth and the wave amplitude respectively.

The SGN system with the mild bottom approximation can be written instead as

\begin{subequations}\label{eq:gnb_mildbot}
\begin{align}
& \partial_t h+ \partial_x(h\overline u) = 0, \label{eq:gnb_mass_mildbot} \\
 & \partial_t (h  \overline u) + \partial_x \left( h{ \overline u}^2 +h \overline p \right)+gh\partial_x h  +\left(gh + \frac{3}{2}  \overline p\right)\partial_x z_b=0, \label{eq:gnb_momx_mildbot}\\
& \partial_t (h \overline w) + \partial_x(h \overline u \, \overline w) = \frac{3}{2}  \overline p, \label{eq:gnb_momz_moldbot}\\
& \overline w + \frac{1}{2}h\partial_x  \overline u -  \overline u\partial_x z_b = 0,\label{eq:p_mildbot}
\end{align} 
\end{subequations}
and its hyperbolic reformulation, as proposed in \cite{escalante:2018,EM20}, reads 

\begin{subequations}\label{eq:hgnb_mildbot}
\begin{align}
& \partial_t h+ \partial_x(h\overline u) = 0, \label{eq:hgnb_mass_mildbot} \\
 & \partial_t (h  \overline u) + \partial_x \left( h{ \overline u}^2 +h \overline p \right)+gh\partial_x h  +\left(gh + \frac{3}{2}  \overline p\right)\partial_x z_b=0, \label{eq:hgnb_momx_mildbot}\\
& \partial_t (h \overline w) + \partial_x(h \overline u \, \overline w) = \frac{3}{2}  \overline p, \label{eq:hgnb_momz_moldbot}\\
& \partial_t\left(h\overline p\right) + \partial_x\left(h\overline u \overline p \right) +c^2\left(h\partial_x \overline u +2\left(\overline w-\overline u\partial_x z_b\right)\right) = 0.\label{eq:hp_mildbot}
\end{align} 
\end{subequations}
Here, $c$ is an \textit{artificial sound speed}, and for $c^2 \to \infty$ it is obvious that \eqref{eq:hp_mildbot}
reduces to the original equation \eqref{eq:p_mildbot}. The underlying idea of the hyperbolic reformulation could be seen as an extension of the \textit{artificial compressibility method} \cite{chorin1} to non-hydrostatic shallow water flows.   
Similarly to what has been done in \cite{escalante:2018} for the SGN system with the mild bottom approximation (\ref{eq:gnb_mildbot}) and for the system derived in \cite{saintmarie:2015}, we propose the following hyperbolic reformulation of the SGN system (\ref{eq:gnb2}).  
We first rewrite system (\ref{eq:gnb2}) as

\begin{subequations}\label{eq:gnb3}
\begin{align}
& \partial_t h+ \partial_x(h\overline u) = 0, \label{eq:gnb_mass3} \\
& \partial_t (h\overline u) + \partial_x \left( h\overline u^2 +h\overline p \right)+gh\partial_x h  + (gh + p_b)\partial_x z_b =0, \label{eq:gnb_momx3}\\
& \partial_t (h\overline w) + \partial_x(h \overline u\overline w) = p_b, \label{eq:gnb_momz3}\\
& \partial_t(h\sigma) + \partial_x(h\overline u\sigma)= -6 p_b +12 \overline p  \label{eq:gnb_43}, \\
& \sigma = -h \partial_x \overline u \label{eq:gnb_53} \\
& \overline w - \frac{\sigma}{2} - \overline u\partial_x z_b = 0, \label{eq:gnb_63} 
\end{align} 
\end{subequations}
where the only difference with respect to (\ref{eq:gnb2}) is that, in the last equation, the term $h\partial_x \overline u /2$ has been rewritten as $-\sigma/2$. We then multiply equations \eqref{eq:gnb_53} and \eqref{eq:gnb_63} by $c^2$ and by $6c^2,$ respectively, and add the corresponding time derivatives of the non-hydrostatic pressure variables so as to obtain two evolution equations for $\overline p$ and $p_b$. The new hyperbolic reformulation of the Serre-Green-Naghdi model proposed in this paper then reads  

\begin{subequations}\label{eq:hgnb}
\begin{align}
& \partial_t h + \partial_x(h\overline u) = 0, \label{eq:hgnb_mass} \\
& \partial_t(h\overline u)      +\partial_x\left(h \overline u^2+h\overline p\right)+gh\partial_x h +(gh+p_b)\partial_x z_b = 0, \label{eq:hgnb_momx} \\
& \partial_t(h\overline w) + \partial_x(h\overline u \ \overline w) =p_b,  \label{eq:hgnb_momz}	\\
& \partial_t(h\sigma) + \partial_x(h \overline u \sigma)  = -6p_b+12\overline p, \label{eq:hgnb_sigma}\\
& \partial_t(h\overline p)+\partial_x[h\overline u(\overline p+c^2)] - c^2 \overline u \partial_x h  = -c^2 \sigma, \label{eq:hgnb_p}\\
& \partial_t(hp_b)+\partial_x (h\overline u p_b )  -6c^2 \overline u \partial_x z_b  = -6c^2 \left(\overline w - \frac{\sigma}{2} \right). \label{eq:hgnb_pb}
\end{align}
\end{subequations}
It is important to notice that for $c^2 \to \infty$ the system \eqref{eq:hgnb} reduces to system \eqref{eq:gnb3}.
Moreover, as it is shown in sections \ref{sec:energybalance} and \ref{sec:eigen}, system \eqref{eq:hgnb} fulfils two other important properties. First, and contrary to the SGN system with the mild bottom approximation and the corresponding hyperbolic reformulation proposed in \cite{escalante:2018}, an energy conservation law can be associated to it. Moreover, it is strictly hyperbolic, which means that it can be efficiently solved employing classical high order discontinuous Galerkin or finite volume methods for hyperbolic PDE. 

\subsection{Energy balance}
\label{sec:energybalance}
We consider the original (non-hyperbolic) SGN system \eqref{eq:gnb2} and we first rewrite equations \eqref{eq:gnb_momx2}, \eqref{eq:gnb_momz2} and \eqref{eq:gnb_52} using the mass conservation equation \eqref{eq:gnb_mass2} as
\begin{subequations}\label{eq:gnb2_rewrite}
\begin{align}
& h\partial_t \overline u + h\overline u \partial_x \overline u+ \partial_x \left( h\overline p\right) +gh\partial_x h + (gh + p_b)\partial_x z_b =0, \label{eq:gnb_momx2_rewrite}\\
& h\partial_t \overline w + h\overline u \partial_x \overline w = p_b, \label{eq:gnb_momz2_rewrite}\\
& h\partial_t \sigma + h\overline u \partial_x\sigma= -6 p_b +12 \overline p
\label{eq:gnb_52_rewrite}. 
\end{align} 
\end{subequations}
We then compute the quantity
\begin{equation}
\overline u (\ref{eq:gnb_momx2_rewrite})+ \overline w(\ref{eq:gnb_momz2_rewrite})+\frac{1}{12}\sigma(\ref{eq:gnb_52_rewrite}) = 0
\label{eq:energy_computation}
\end{equation}
and using again the mass conservation equation \eqref{eq:gnb_mass2}, yields 
\begin{multline}
\partial_t\left[\frac{h}{2} \left( \overline u^2 +\overline w^2 +\frac{\sigma^2}{12}\right)  \right] + \partial_x\left[\overline u\frac{h}{2}\left(\overline u^2 + \overline w^2 + \frac{\sigma^2}{12}\right) + \overline u\frac{gh^2}{2} +\overline u h \overline p \right] \\   -h\overline p \partial_x \overline u -\frac{g h^2}{2} \partial_x \overline u + \overline u\left(gh+p_b\right)\partial_x z_b -\overline w p_b + \frac{\sigma p_b}{2} -\sigma \overline p = 0.
\label{eq:both_systems}
\end{multline}
We now use equations  \eqref{eq:gnb_53} and \eqref{eq:gnb_63} for $\overline w$ and $\sigma$, obtaining
\begin{multline*}
\partial_t\left[\frac{h}{2} \left( \overline u^2 +\overline w^2 +\frac{\sigma^2}{12}\right)  \right] + \partial_x\left[\overline u\frac{h}{2}\left(\overline u^2 + \overline w^2 + \frac{\sigma^2}{12}\right) + \overline u\frac{gh^2}{2} +\overline u h \overline p\right]   -\frac{gh^2}{2}\partial_x \overline u + gh\overline u \partial_x z_b = 0.
\end{multline*}
Adding the quantity $gh(\partial_t h +\partial_x(hu))/2$ , which is equal to zero thanks to equation (\ref{eq:gnb_mass2}), and using $\partial_t z_b = 0$ we get
\begin{equation*}
\partial_t\left(\frac{h}{2} \left[ \overline u^2 +\overline w^2 +\frac{\sigma^2}{12}\right] + g\left(h+2z_b\right)  \right) + \partial_x\left[\overline u\frac{h}{2}\left(\overline u^2 + \overline w^2 + \frac{\sigma^2}{12}+g\left(h+2z_b\right)\right)+ \overline u\frac{gh^2}{2} +\overline u h \overline p\right]  = 0.
\end{equation*}
If we introduce the following definition for the energy:
\begin{equation}
\tilde{E} = \frac{h}{2}\left[\overline u^2+\overline w^2+\frac{\sigma^2}{12}+g(h+2z_b) \right],
\label{eq:energy_gnb}
\end{equation}
we obtain the energy conservation equation
\begin{equation}
\partial_t \tilde{E} + \partial_x\left[\overline u\left(\tilde{E}+\frac{gh^2}{2} + h\overline p\right)\right] = 0.
\end{equation}

We now consider the hyperbolic system \eqref{eq:hgnb}, for which equation \eqref{eq:both_systems} is also valid. 
We then rewrite equations \eqref{eq:hgnb_p} and \eqref{eq:hgnb_pb} as follows:
\begin{subequations}
\begin{align}
& \overline p\left(-h\partial_x \overline u -\sigma \right) = \frac{1}{c^2}\left(\partial_t \left(\frac{h\overline p^2}{2}\right) + \partial_x\left(\frac{h\overline u\, \overline p^2}{2}\right)\right)\\
& p_b\left(\overline u \partial_x z_b -\overline w + \frac{\sigma}{2}\right) = \frac{1}{6c^2}\left(\partial_t \left(\frac{h p_b^2}{2}\right) + \partial_x\left(\frac{h\overline u  p_b^2}{2}\right)\right),
\end{align}
\end{subequations}
and we substitute them into equation \eqref{eq:both_systems}, obtaining
\begin{multline}
\partial_t\left[\frac{h}{2} \left( \overline u^2 +\overline w^2 +\frac{\sigma^2}{12} + \frac{\overline p^2}{c^2}+\frac{p_b^2}{6c^2}\right)  \right] \\+ \partial_x\left[\overline u\frac{h}{2}\left(\overline u^2 + \overline w^2 + \frac{\sigma^2}{12} + \frac{\overline p^2}{c^2}+\frac{p_b^2}{6c^2}\right) + \overline u\frac{gh^2}{2} +\overline u h \overline p\right]   -\frac{gh^2}{2}\partial_x \overline u + gh\overline u \partial_x z_b = 0.
\end{multline}
As for the non-hyperbolic system, we add the quantity $gh(\partial_t h +\partial_x(hu))/2$ and we use $\partial_t z_b = 0$, thus retrieving the energy balance 

\begin{equation}
\partial_t E + \partial_x\left(\overline u E+\frac{gh^2}{2} + h\overline p\right) = 0,
\label{eq:energy_hgnb_law}
\end{equation}
where $E$ is now defined as
\begin{equation}
E = \frac{h}{2}\left[\overline u^2+\overline w^2+g(h+2z_b) +\frac{\sigma^2}{12}+ \frac{\overline p^2}{c^2}+ \frac{p_b^2}{6c^2}\right].
\label{eq:energy_hgnb}
\end{equation}
If we compare equations (\ref{eq:energy_gnb}) and (\ref{eq:energy_hgnb}), we can notice that the energy  $E$ associated to the hyperbolic system reduces to the energy $\tilde{E}$  associated to the original system for $c^2 \to \infty$.

\subsection{Eigenstructure of the hyperbolic reformulation of the Serre-Green-Naghdi system}
\label{sec:eigen}
Defining the vector of unknowns as   
\begin{equation}
\mathbf{U} = (h, h\overline u, h\overline w, h\sigma, h\overline p, h p_b)
\end{equation}
we rewrite the hyperbolic SGN system (\ref{eq:hgnb}) in compact form as follows
\begin{equation}
\partial_t \mathbf{U} + \partial_x\mathbf{F}(\mathbf{U}) + \mathbf{B}(\mathbf{U}) \partial_x\mathbf{U} = \mathbf{S}(\mathbf{U}),
\label{eq:general_pde}
\end{equation}
where $\mathbf{F}(\mathbf{U})$ denotes the nonlinear flux tensor, $\boldsymbol{\mathbf{B}}(\mathbf{U}) \cdot\partial_{x}
\mathbf{U}$ is a genuinely non-conservative term and $\mathbf{S}\left( \mathbf{U}\right)$ corresponds to the source term,
\begin{align*}
& \mathbf{F}(\mathbf{U}) = (h\overline u, h\overline u^2 + h \overline p , h \overline u\,\overline w, h \overline u \sigma, h \overline u (\overline p + c^2), h \overline u p_b), \nonumber\\
&\mathbf{B}(\mathbf{U})\partial_x \mathbf{U}= (0, gh\partial_x h+(gh + p_b)\partial_x z_b, 0, 0, -c^2 \overline u \partial_x h, -6c^2\overline u \partial_x z_b ), \nonumber	\\
& \mathbf{S}(\mathbf{U})= (0, 0, p_b, -6p_b + 12\overline p, -c^2 \sigma, -6c^2 (\overline w-\sigma/2)). \nonumber
\end{align*}
\label{eq:fluxes}
System (\ref{eq:general_pde}) can be rewritten in quasilinear form as
\begin{equation}
\partial_t \mathbf{U} + \mathbf{A}(\mathbf{U})\partial_x \mathbf{U} = \mathbf{S}(\mathbf{U}), 
\label{eqn.ql1d} 
\end{equation}
with $ \mathbf{A}(\mathbf{U}) = \mathbf{J}_{\mathbf{F}}+ \mathbf{B}(\mathbf{U})$ and $\mathbf{J}_{\mathbf{F}} = \partial \mathbf{F}/\partial \mathbf{U}$.
Matrix $\mathbf{A}(\mathbf{U})$ has six real eigenvalues $\lambda_{1,2,3,4}=u$, $\lambda_{5,6}= u \pm  \sqrt{\overline p +gh +c^2}$, while the corresponding set of linearly independent eigenvectors is
\begin{align*}
& \mathbf{r}_1 = (0, 0, 0, 1, 0, 0), \\
& \mathbf{r}_2 = (0, 0, 0, 0, 0, 1), \\
& \mathbf{r}_3 = (0, 0, 1, 0, 0, 0), \\
& \mathbf{r}_4 = (1, \overline u, 0, 0, -gh, 0), \\
& \mathbf{r}_5 = (1, \overline u + \sqrt{\overline p +gh +c^2}, \overline w, \sigma, \overline p +c^2, p_b), \\
& \mathbf{r}_6 = (1, \overline u - \sqrt{\overline p +gh +c^2}, \overline w, \sigma, \overline p +c^2, p_b).
\end{align*}
Therefore, the system \eqref{eq:hgnb} is hyperbolic and satisfies the extra energy conservation law 
\eqref{eq:energy_hgnb_law}. 

\subsection{Two-dimensional extension of the model}
The two-dimensional model is simply obtained by adding a second depth-averaged velocity component $\bar{v}$ to the system. The full set of governing equations in the two-dimensional case therefore reads 
\begin{equation}
   \partial_t \mathbf{U} + \nabla \cdot \mathbf{F}(\mathbf{U}) + \mathbf{B}(\mathbf{U}) \cdot \nabla \mathbf{U} = \mathbf{S}(\mathbf{U}), 
\end{equation}
with the state vector 
\begin{equation}
   \mathbf{U} = (h, h\overline u, h\overline v, h\overline w, h\sigma, h\overline p, h p_b), 
\end{equation} 
the flux tensor  $ \mathbf{F}(\mathbf{U}) = (  \mathbf{F}_1 ,  \mathbf{F}_2)$ with 
\begin{eqnarray}
  \mathbf{F}_1 & =  (h\overline u, h\overline u^2 + h \overline p , h \overline u \, \overline v, h \overline u\,\overline w, h \overline u \sigma, h \overline u (\overline p + c^2), h \overline u p_b), & \\ 
  \mathbf{F}_2 & = (h\overline v, h \overline v \, \overline u, h\overline v^2 + h \overline p, 
  h \overline v\,\overline w, h \overline v \sigma, h \overline v (\overline p + c^2), h \overline v p_b), & 
\end{eqnarray}
the non-conservative product $\mathbf{B}(\mathbf{U}) \cdot \nabla \mathbf{U} = \mathbf{B}_1(\mathbf{U}) \partial_x \mathbf{U} + \mathbf{B}_2(\mathbf{U}) \partial_y \mathbf{U} $ with  
\begin{eqnarray}
\mathbf{B}_1(\mathbf{U}) \partial_x \mathbf{U} & = & 
(0, gh\partial_x h+(gh + p_b)\partial_x z_b, 0, 0, 0, -c^2 \overline u \partial_x h, -6c^2\overline u \partial_x z_b ),  \\ 
\mathbf{B}_2(\mathbf{U}) \partial_y \mathbf{U} & = & (0, 0, gh\partial_y h+(gh + p_b)\partial_y z_b, 0, 0, -c^2 \overline v \partial_y h, -6c^2\overline v \partial_y z_b ),     
\end{eqnarray}
and the algebraic source term 
\begin{equation}
    \mathbf{S}(\mathbf{U})= (0, 0, 0, p_b, -6p_b + 12\overline p, -c^2 \sigma, -6c^2 (\overline w-\sigma/2)).
\end{equation}
It is easy to check that in the multi-dimensional case, the extra energy conservation law reads 
\begin{equation}
\partial_t E + \partial_x\left(\overline u E+\frac{gh^2}{2} + h\overline p\right) 
             + \partial_y\left(\overline v E+\frac{gh^2}{2} + h\overline p\right) = 0, 
\label{eq:energy_2d}
\end{equation}
with the energy
\begin{equation}
E = \frac{h}{2}\left( \overline u^2 + \overline v^2 + \overline w^2+g(h+2z_b) +\frac{\sigma^2}{12}+ \frac{\overline p^2}{c^2}+ \frac{p_b^2}{6c^2}\right).
\label{eq:energy_2d}
\end{equation}

\section{Numerical scheme} \label{sec:numscheme}
To study the behaviour of the solutions of the new model \eqref{eq:hgnb} and compare it with the solutions 
provided by the hyperbolic reformulation of the SGN system with mild bottom assumption \cite{escalante:2018}, 
we employ the high order accurate fully-discrete one-step ADER discontinuous Galerkin methodology. 
The family of ADER-DG schemes has been developed during the last decades for both unstructured and Cartesian mesh including space-time meshes, see   \cite{ADERNC,ADERNSE,DumbserZanotti,Dumbser2014,Zanotti2015,GPRmodel,DFTBW18,BCDGP20}. 
The results obtained so far for a wide variety of hyperbolic models, from the compressible Navier-Stokes equations to general relativity and the unified Godunov-Peshkov-Romenski 
model of continuum mechanics, make ADER-DG a suitable candidate to be used also for the discretization of the 
new hyperbolic model proposed in this paper that concerns non-hydrostatic flows over general bottom topographies. In what follows, we provide a brief description of the method in 2D, for further details we refer 
to the above references. 

\subsection{Fully discrete one-step ADER-DG schemes}

We consider the computational domain, $\Omega$, to be covered using a Cartesian grid 
whose elements are of the form $\Omega_{i}=\left[x_{i}-\frac{1}{2}\Delta x,x_{i}+\frac{1}{2}\Delta x\right]\times \left[y_{i}-\frac{1}{2}\Delta y,y_{i}+\frac{1}{2}\Delta y\right]$ 
with $\mathbf{x}_{i}=\left(x_{i},y_{i}\right) $ the barycentre of cell $\Omega_{i}$ 
and $\Delta x$, $\Delta y$ the cell size in each spatial coordinate direction.
Denoting by $\mathbf{u}_h(\mathbf{x},t^n)$ the discrete solution of \eqref{eq:general_pde} 
written in the space of piecewise polynomials of degree $N$, the discrete solution is sought under 
the form
\begin{equation}
\mathbf{u}_h(\mathbf{x},t^n) = \phi_l (\mathbf{x})\;
\hat{\mathbf{u}}^n_l,   \quad \mathbf{x} \in
\Omega_i. \label{eq:disc_sol}
\end{equation}
Here, $\phi_l(\mathbf{x}) = \phi_{l_1}(\xi) \phi_{l_2}(\eta)$ are the basis functions, which are 
chosen to be tensor products of one-dimensional basis functions $\phi_{l_m}(\chi)$ on the unit
interval $\chi \in \Omega_{\mathrm{ref}}=\left[0,1\right]$. The reference coordinates 
$0 \leq \xi, \eta \leq 1$ are obtained via the transformations 
 $x = x_{i}-\frac{1}{2}\Delta x + \xi \Delta x$ and 
 $y = y_{i}-\frac{1}{2}\Delta y + \eta \Delta y$, respectively. Throughout this paper 
the classical Einstein summation convection is used and $l$ is a multidimensional 
index, referring to the one-dimensional basis functions $\phi_{l_{m}}$ 
to be used in the tensor product.  
In particular, we consider the Lagrange interpolation polynomials passing through the 
Gauss-Legendre quadrature points of a $N+1$ Gaussian quadrature formula, which are by construction 
orthogonal.  
As a consequence of the nodal tensor-product basis employed here, the scheme can be written in a 
dimension by  dimension fashion and integral operators are decomposed in the product of 
one-dimensional operators. 

We now multiply the governing PDE system \eqref{eq:general_pde}, by test functions $\phi_k$, which are 
identical to the basis functions, and integrate it over a space-time control volume $[t^{n},t^{n+1}]\times\Omega_i$, obtaining the following weak problem
\begin{equation}
\int \limits_{t^n}^{t^{n+1}} \int\limits_{\Omega_i } \phi_k \left(
\partial_t \mathbf{U}  + \nabla \cdot
\mathbf{F}(\mathbf{U})+ \boldsymbol{\mathbf{B}}(\mathbf{U}) \cdot\nabla
\mathbf{U}\right) \,d\mathbf{x}\,dt = 
\int \limits_{t^n}^{t^{n+1}} \int\limits_{\Omega_i } \phi_k \, \mathbf{S}\left( \mathbf{U} \right).\label{eq:weakPDE}
\end{equation}
By taking into account \eqref{eq:disc_sol}, integrating the flux divergence term by 
parts in space and the time derivative by parts in time the above weak problem becomes 
\begin{gather}
\left( \, \int\limits_{\Omega_i} \phi_k \phi_l \, d\mathbf{x}\right)
\left( \hat{\mathbf{u}}_l^{n+1} - \hat{\mathbf{u}}_l^{n} \, \right) +
\int \limits_{t^n}^{t^{n+1}} \int\limits_{\partial
	\Omega_i } 
\phi_k \mathcal{G}\left(\q_h^-, \q_h^+ \right) \cdot
\mathbf{n} \, dS \, dt 
+ \int \limits_{t^n}^{t^{n+1}} \!\!
\int\limits_{\partial \Omega_i } \!\! \phi_k \mathcal{D}\left(
\q_h^-,\q_h^+ \right) \cdot \mathbf{n} \, dS \, dt \nonumber \\  -
\int \limits_{t^n}^{t^{n+1}} \!\! \int\limits_{\Omega_i } \!\!\! \nabla
\phi_k \cdot \mathbf{F}(\q_h) \,d\mathbf{x}\,dt + \int
\limits_{t^n}^{t^{n+1}} \!\! \int\limits_{\Omega_i^{\circ} } \phi_k
\mathbf{\mathbf{B}}(\mathbf{q}_h) \cdot \nabla \mathbf{q}_h
\,d\mathbf{x}\,dt = 
\int\limits_{\Omega_i  } \phi_k
\mathbf{\mathbf{S}}(\mathbf{q}_h)  
\,d\mathbf{x}\,dt,  
\label{eq:ADER-DG}
\end{gather}
where we have denoted by $\mathbf{n}$ the outward unit normal at the cell boundary 
$\partial \Omega_{i}$, and $\q_h$ is a local space-time predictor whose computation 
will be detailed in the next section. Let us remark that the test and basis functions 
can jump across the element interfaces, leading to the so-called discontinuous Galerkin 
finite element method. To account for these jumps, we make use of Riemann solvers at 
the element interfaces,  see e.g. \cite{toro-book} for a broad overview of different 
exact and approximate Riemann solvers. In this paper, we use either the simple Rusanov-type 
flux  
\begin{equation}
\mathcal{G}\left(\q_h^-, \q_h^+ \right) \cdot \mathbf{n} = \frac{1}{2}
\left( \mathbf{F}(\q_h^+) + \mathbf{F}(\q_h^-) \right) \cdot \mathbf{n}
- \frac{1}{2} s_{\max} \left( \q_h^+ - \q_h^- \right)\,,
\label{eq.rusanov} 
\end{equation} 
with the maximum wavespeed at the interface $s_{\max}$, 
or the more sophisticated generalized Osher-type scheme forwarded in \cite{OsherNC,OsherUniversal}.  
Here, $\q_h^-$ and $\q_h^+$ denote the boundary-extrapolated values of the predictor from within the
element and from the neighbor element, respectively. 
Furthermore, the scheme also requires a proper  
discretization of the non conservative products at the boundaries arising due to the 
presence of the bottom slope term. To this end, we consider the works on so-called path 
conservative schemes forwarded by Castro, Par\'es and collaborators in \cite{Castro2006,Pares2006,Castro2d,Munoz2007,MunozPares,Castro2008}, which are  
based on the theory of Dal Maso, Le Floch and Murat \cite{DLMtheory} on nonconservative 
hyperbolic PDE systems. For a more detailed discussion on the topic, see also \cite{NCproblems}
and references therein. 
The first extensions of the path-conservative approach to higher order DG schemes can be found in \cite{Rhebergen2008,ADERNC}.   
Within the path-conservative framework, it is very simple and natural to construct also so-called 
well-balanced schemes for shallow water models, see \cite{Bermudez1994,GarciaNavarro1,LeVequeWB}. 
The jump term in the non-conservative product is computed using a path integral in space between the two 
extrapolated values related to the face, $\q_h^-$ and $\q_h^+$,
\begin{equation}
\mathcal{D}\left( \q_h^-,\q_h^+ \right) \cdot \mathbf{n} =  \frac{1}{2}
\left(\int \limits_{0}^{1} \mathbf{\mathcal{B}} \left(
\bm{\psi}(\q_h^-,\q_h^+,s) \right)\cdot\mathbf{n} \, ds
\right)\cdot\left(\q_h^+ - \q_h^-\right)\,, \label{eq:PC}
\end{equation}
where we use the linear segment path
\begin{equation}
\bm{\psi} = \bm{\psi}(\q_h^-, \q_h^+, s) = \q_h^- + s \left( \q_h^+ -
\q_h^- \right)\,, \qquad s \in [0,1].
\end{equation}

\subsection{Local space-time predictor}
To determine the local space-time predictor solution, $\q_h(\mathbf{x},t)$, which will 
lead to a high order scheme in space and time avoiding the cumbersome Cauchy-Kovalewskaya 
procedure used in original ADER schemes \cite{Toro2001,Toro2002,toro3,TT04}, we employ the weak formulation 
in space-time proposed in \cite{DumbserEnauxToro,Dumbser2008}. Consequently, the Cauchy problem 
is solved ``in the small'' thus neglecting the iteration between neighbours. 

Let us consider a space-time test function, $\theta_{k}=\theta_{k}(\mathbf{x},t) $, 
built as the product of the one dimensional spatial basis functions already introduced 
and an additional nodal basis function for the time dependency. Multiplying 
\eqref{eq:general_pde} by this test function and integrating over the space-time 
control volume, $\Omega_{i}\times\left[t^{n},t^{n+1}\right]$, we get
\begin{equation}
\int \limits_{t^n}^{t^{n+1}} \!\! \int\limits_{\Omega_i  }
\theta_k \, \partial_t \q_h \,d\mathbf{x} \, dt + \int
\limits_{t^n}^{t^{n+1}} \!\! \int\limits_{\Omega_i  } \theta_k \,
\nabla \cdot \mathbf{F}(\q_h) \,d\mathbf{x}\,dt  
+ \int
\limits_{t^n}^{t^{n+1}} \!\! \int\limits_{\Omega_i^{\circ} } \theta_k
\mathbf{\mathbf{B}}(\q_h ) \cdot \nabla \q_h \,d\mathbf{x}\,dt =
\int\limits_{\Omega_i } \theta_k
\mathbf{S}(\q_h )   \,d\mathbf{x}\,dt. 
\label{eq:predictor1}
\end{equation}
Integration by parts in time just of the first term yields 
\begin{eqnarray}
\int\limits_{\Omega_i }
\theta_k(\mathbf{x},t^{n+1}) \q_h(\mathbf{x},t^{n+1}) \,d\mathbf{x}  -  
\int\limits_{\Omega_i }
\theta_k(\mathbf{x},t^{n}) \mathbf{u}_h(\mathbf{x},t^{n}) \,d\mathbf{x}     
-\int \limits_{t^n}^{t^{n+1}} \!\! \int\limits_{\Omega_i  }
\partial_t \theta_k \,  \q_h \,d\mathbf{x} \, dt + 
\nonumber \\ 
\int
\limits_{t^n}^{t^{n+1}} \!\! \int\limits_{\Omega_i } \theta_k \,
\nabla \cdot \mathbf{F}(\q_h) \,d\mathbf{x}\,dt  
+ \int
\limits_{t^n}^{t^{n+1}} \!\! \int\limits_{\Omega_i^{\circ} } \theta_k
\mathbf{\mathbf{B}}(\q_h ) \cdot \nabla \q_h \,d\mathbf{x}\,dt =
\int\limits_{\Omega_i } \theta_k
\mathbf{S}(\q_h )   \,d\mathbf{x}\,dt, 
\label{eq:predictor}
\end{eqnarray}
that corresponds to a nonlinear system from which the unknown 
degrees of freedom $ {\hat{\q}}_k$ of the space-time expansion, 
\begin{equation}
\q_h(\boldsymbol{x},t) = \theta_k(\boldsymbol{x},t) \, {\hat{\q}}_k,
\label{eq.stdof}
\end{equation}
can be computed after substitution of the value obtained for each spatial
degree of freedom at the previous time step, $\mathbf{u}_h(\mathbf{x},t^{n})$. 
In \eqref{eq.stdof} the $\theta_k(\boldsymbol{x},t) = \phi_{l_0}(\tau) \phi_{l_1}(\xi) \phi_{l_2}(\eta)$ 
are nodal space-time basis functions, which are again taken to be tensor products of the 
one-dimensional Lagrange interpolation polynomials $\phi_{l_m}(\chi)$,  
passing through the Gauss-Legendre quadrature points on the unit interval, with the additional 
transformation for the reference time $\tau$ given by $t = t^n + \tau \Delta t$.   
The solution of \eqref{eq:predictor} can be found via a fast-converging iterative
fixed point scheme, the convergence of which was proven in \cite{BCDGP20}.  

\section{Numerical tests} \label{sec:numtest}
In this section, the new model (\ref{eq:hgnb}) proposed in this paper is tested at the aid of several numerical experiments, including comparisons with quasi exact, numerical and experimental reference solutions. 
The first considered  benchmark analyses the propagation of a solitary wave over a flat bottom. This test aims at validating both the mathematical model and the numerical scheme: for this reason also a numerical convergence study has been carried out. 
After having verified the correct behaviour of both the numerical scheme and the mathematical model, different test cases with non trivial bottom topography have also been considered. We have focused especially
on strongly varying topographies,  since the peculiarity of system \eqref{eq:hgnb} is that it is derived from the SGN system for general bottom, unlike the models presented in \cite{JSMarie,escalante:2018}. 
In the following, all the results are presented employing the international system of units (SI).

\subsection{Solitary wave over a flat bottom}
The test case concerning  a solitary wave over a flat bottom  is defined on the computational domain $\Omega = [-50,50]$ with final simulation time $t_{\mathrm{end}}=2$. Concerning boundary conditions, periodic boundaries have been imposed. The initial amplitude of the soliton is $A=0.2$, the still water depth is $H=1$ and, at $t=0$, the soliton is centred in $x_0=0$. The artificial sound velocity $c$ has been set equal to $20$. 

Notice that, in the flat bottom case, the SGN equations without the mild bottom approximation, \eqref{eq:gnb2}, reduce to the SGN equations with the mild bottom approximation, \eqref{eq:gnb_mildbot}. However, concerning the initial condition, we have not employed the analytical solution of system \eqref{eq:gnb_mildbot} (see, for example, \cite{saintmarie:2015}), since this analytical solution, in the case of flat bottom, is an exact solution of system \eqref{eq:gnb2}, but is not an exact solution of the hyperbolic system \eqref{eq:hgnb}. Its use would lead to a deterioration of the order of convergence, especially for the non-hydrostatic averaged pressure $\overline p$ and the non-hydrostatic pressure at the bottom boundary $p_b$. We are instead looking for self similar solutions of the hyperbolic system \eqref{eq:hgnb} of the form  
\begin{equation}
   \mathbf{U}(x,t) = \mathbf{U}(\zeta), \qquad \textnormal{with} \qquad \zeta = x-V \, t, 
\end{equation}
where $V$ is the velocity of the solitary wave and $\zeta$ is the similarity coordinate. Obviously, under this assumption one has $\partial_t \mathbf{U} = -V \, \mathbf{U}'$ and $\partial_x \mathbf{U} = \mathbf{U}'$. Hence, the quasilinear PDE \eqref{eqn.ql1d} can be rewritten as 
\begin{equation}
   - V \, \mathbf{U}' + \mathbf{A}(\mathbf{U}) \mathbf{U}' = \mathbf{S}(\mathbf{U}) 
\end{equation}
and therefore reduces to the following nonlinear ODE system 
\begin{equation}
 \mathbf{U}' = \left( \mathbf{A}(\mathbf{U}) - V \, \mathbf{I} \right)^{-1} \mathbf{S}(\mathbf{U}),  
 \label{eqn.ode.sol} 
\end{equation}
with initial condition $\mathbf{U}(\zeta_0) = ( H_0, 0, 0, 0, \epsilon, 0 )$  
and with $\mathbf{I}$ being the identity matrix. We set $\epsilon = 10^{-8}$ and the ODE system \eqref{eqn.ode.sol} is solved with a $10^{th}$ order DG scheme in time, see \cite{ADERNSE}, in 
order to provide the initial condition for the solitary wave of the new hyperbolic reformulation 
of the SGN system proposed in this paper. For the calculations in this section,  
we use the following parameters: 
$H_0 = 1$, $A=0.2$, $g=9.81$, $c_0 = 20$ and the velocity $V$ is chosen as $V = \sqrt{g (A+H_0)} $. 

In Table \ref{tab:1}, the $L_2$ errors and the convergence rates are shown. Overall, we can observe that the correct order of accuracy is retrieved for the tested polynomial degrees $N=3,5,6,7$, even if, in a few cases, a suboptimal convergence order is obtained.

\begin{table}
\centering
\begin{tabular}{c | c | c c c | c c c | c|}
\toprule
$N$           & $N_x$ & $L_2$ err $h$ & $L_2$ err $\overline u $ & $L_2$ err $\overline p$ & $L_2$ ord $h$ & $L_2$ ord $\overline u$ & $L_2$ ord $\overline p$ & Theor.    \\
\hline
3             &   80 & 6.23E-4 & 4.77E-4 & 5.02E-3 & - & - & - & 4   \\
               & 100 & 2.22E-4 & 1.54E-4 & 2.03E-3 & 4.62 & 5.07 & 4.06 & \\
               & 120 & 9.29E-5 & 5.47E-5 & 1.13E-3 & 4.78 & 5.68 & 4.06 & \\
               & 140 & 4.62E-5 & 2.13E-5 & 4.63E-4 & 4.53 & 6.12 & 4.53 & \\
               & 160 & 2.65E-5 & 9.73E-6 & 2.65E-4 & 4.16 & 5.87 & 4.22 & \\
\hline 
5             &   20 & 5.82E-3 & 8.55E-3 & 5.24E-2 & - & - & - & 6   \\
               &   40 & 2.22E-4 & 1.10E-4 & 1.77E-3 & 4.71 & 6.28 & 4.89 & \\
               &   60 & 2.61E-5 & 1.15E-5 & 2.12E-4 & 5.28 & 5.57 & 5.23 & \\
               &   80 & 4.89E-6 & 2.05E-6 & 3.62E-5 & 5.82 & 5.99 & 6.14 & \\
               & 100 & 1.30E-6 & 4.92E-7 & 1.01E-5 & 5.94 & 6.40 & 5.72 & \\
\hline
6             &   30 & 2.05E-4 & 1.74E-4 & 1.84E-3 & - & - & - & 7   \\
               &   40 & 3.44E-5 & 2.83E-5 & 3.10E-4 & 6.20 & 6.31 & 6.19 & \\
               &   50 & 7.71E-6 & 3.73E-6 & 6.87E-5 & 6.70 & 9.08 & 6.75 & \\
               &   60 & 1.46E-6 & 8.31E-7 & 1.58E-5 & 9.13 & 8.24 & 8.06 & \\
               &   70 & 4.01E-7 & 2.62E-7 & 5.00E-6 & 8.38 & 7.49 & 7.46 & \\
\hline
7             &   10 & 1.92E-2 & 2.98E-2 & 0.11 & - & - & - & 8   \\
               &   20 & 4.99E-4 & 7.17E-4 & 5.77E-3 & 5.27 & 5.38 & 4.25 & \\
               &   30 & 4.09E-5 & 2.99E-5 & 3.57E-4 & 6.17 & 7.84 & 6.86 & \\
               &   40 & 4.62E-6 & 2.15E-6 & 4.44E-5 & 7.58 & 9.15 & 7.25 & \\
               &   50 & 4.31E-7 & 3.44E-7 & 5.61E-6 & 10.63 & 8.21 & 9.27 & \\
\bottomrule
\end{tabular}
\caption{$L_2$ errors and convergence rates for the solitary wave over a flat bottom test case (polynomial degrees $N=3,5,6,7$), at the final simulation time $t_{\mathrm{end}}=2$.}
\label{tab:1}
\end{table}

\begin{figure}
\centering
\begin{subfigure}[$h$]{
     \label{fig:flat_bottom_h}
     \includegraphics[width=0.42\textwidth]{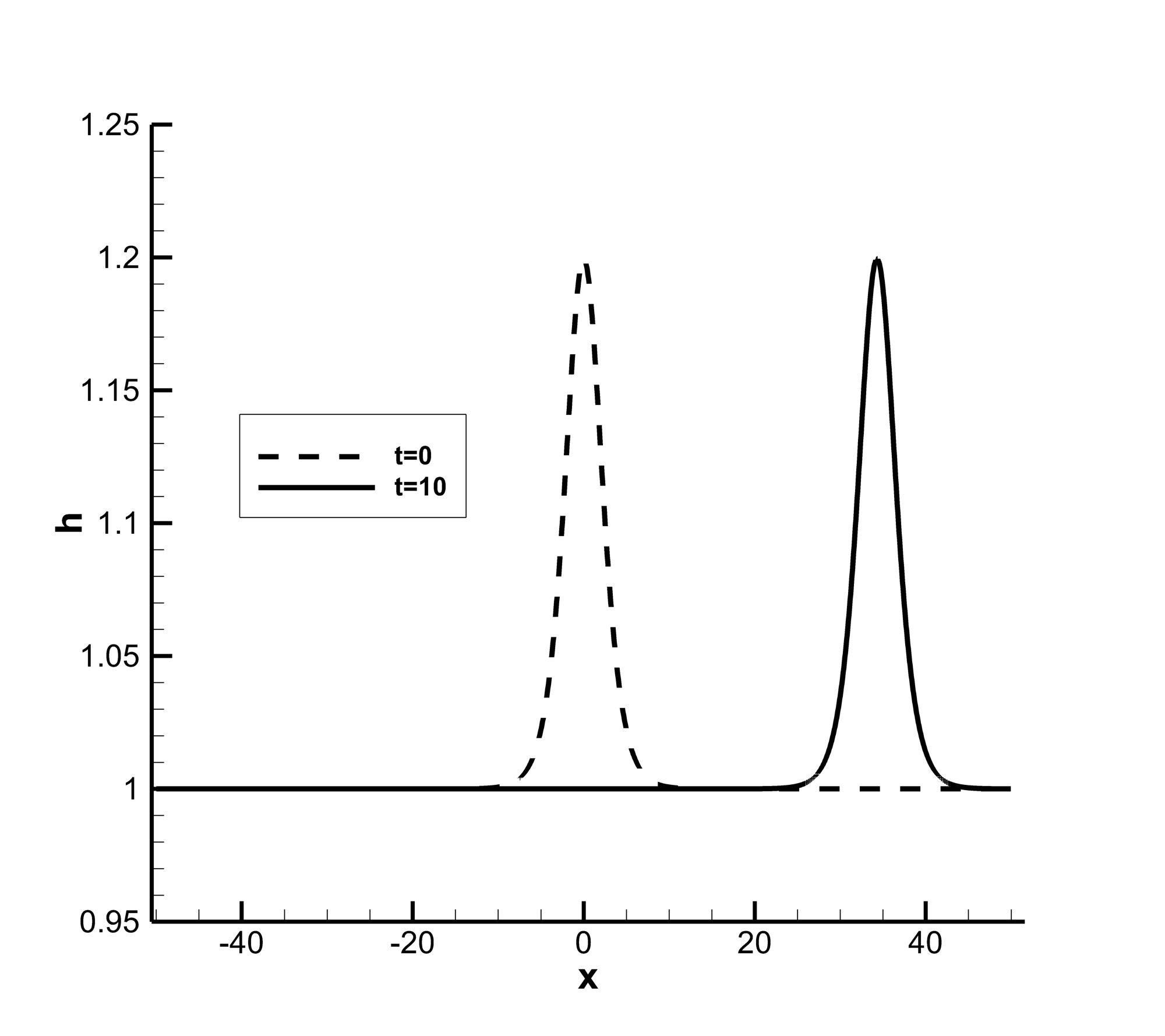}}    
\end{subfigure}
\begin{subfigure}[$\overline u$]{
      \label{fig:flat_bottom_u}
     \includegraphics[width=0.42\textwidth]{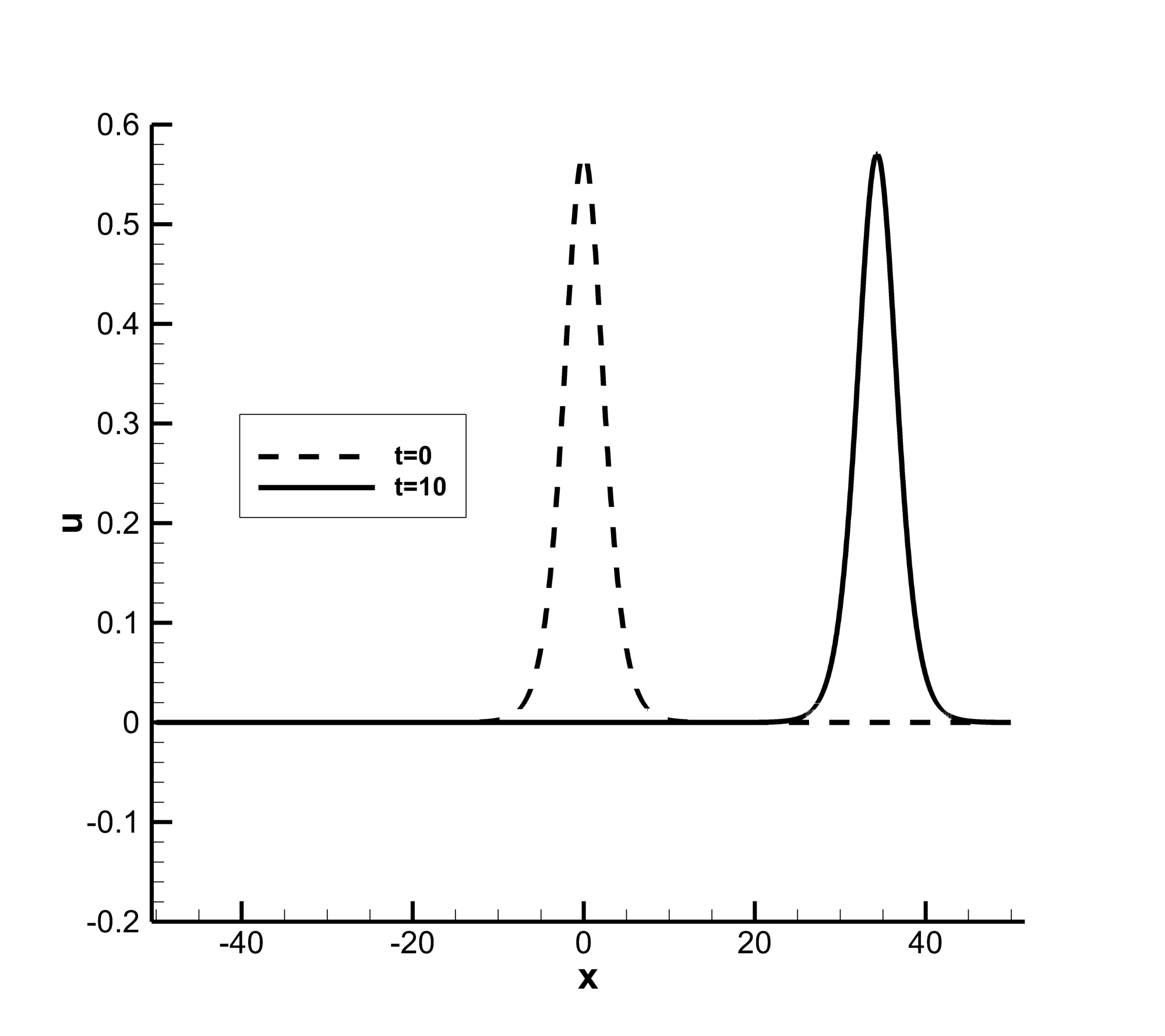}}
\end{subfigure} 

\begin{subfigure}[$\overline w$]{
     \label{fig:flat_bottom_w}
    \includegraphics[width=0.42\textwidth]{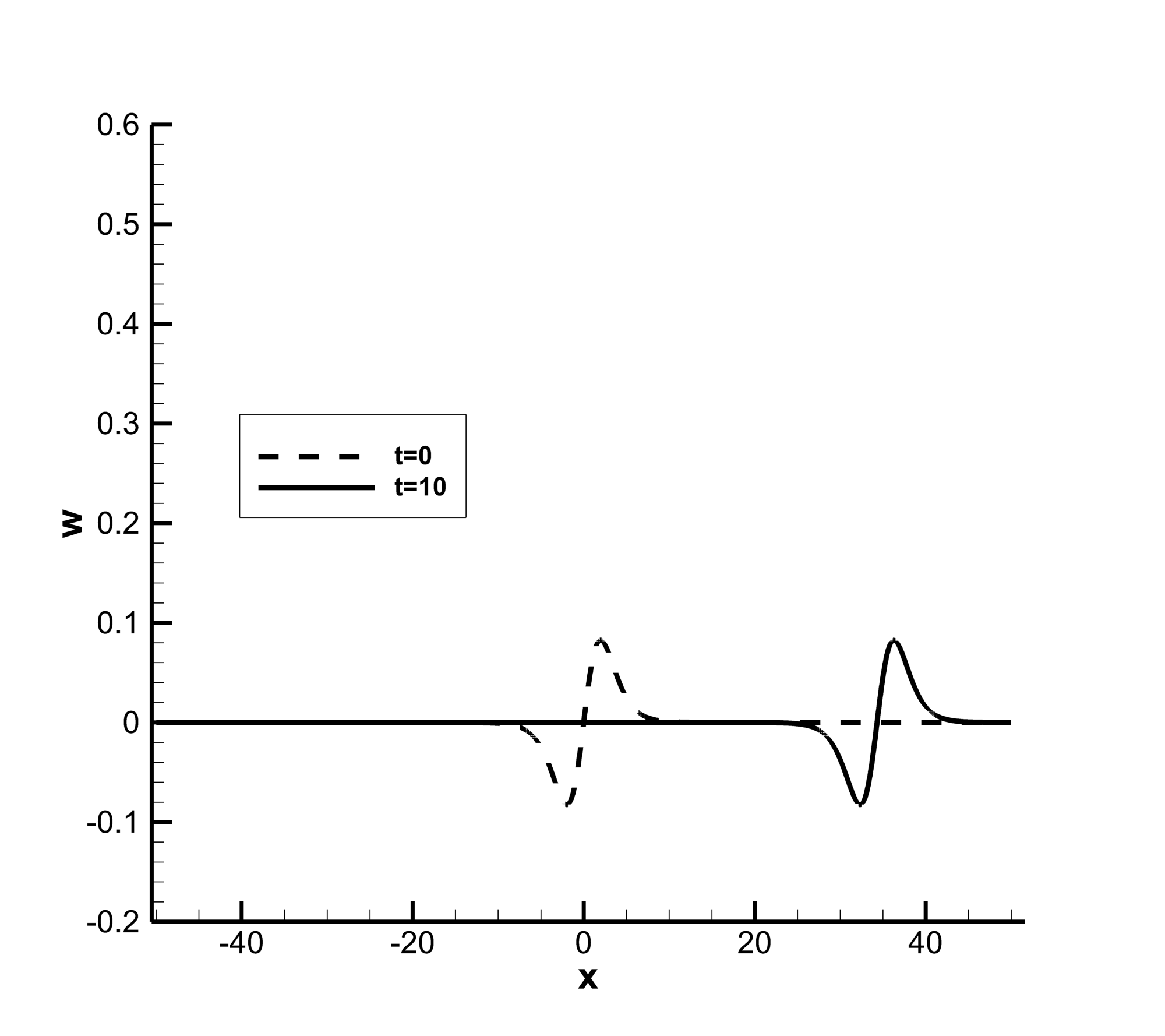}}    
\end{subfigure}
\begin{subfigure}[$\sigma$]{
      \label{fig:flat_bottom_sigma}
     \includegraphics[width=0.42\textwidth]{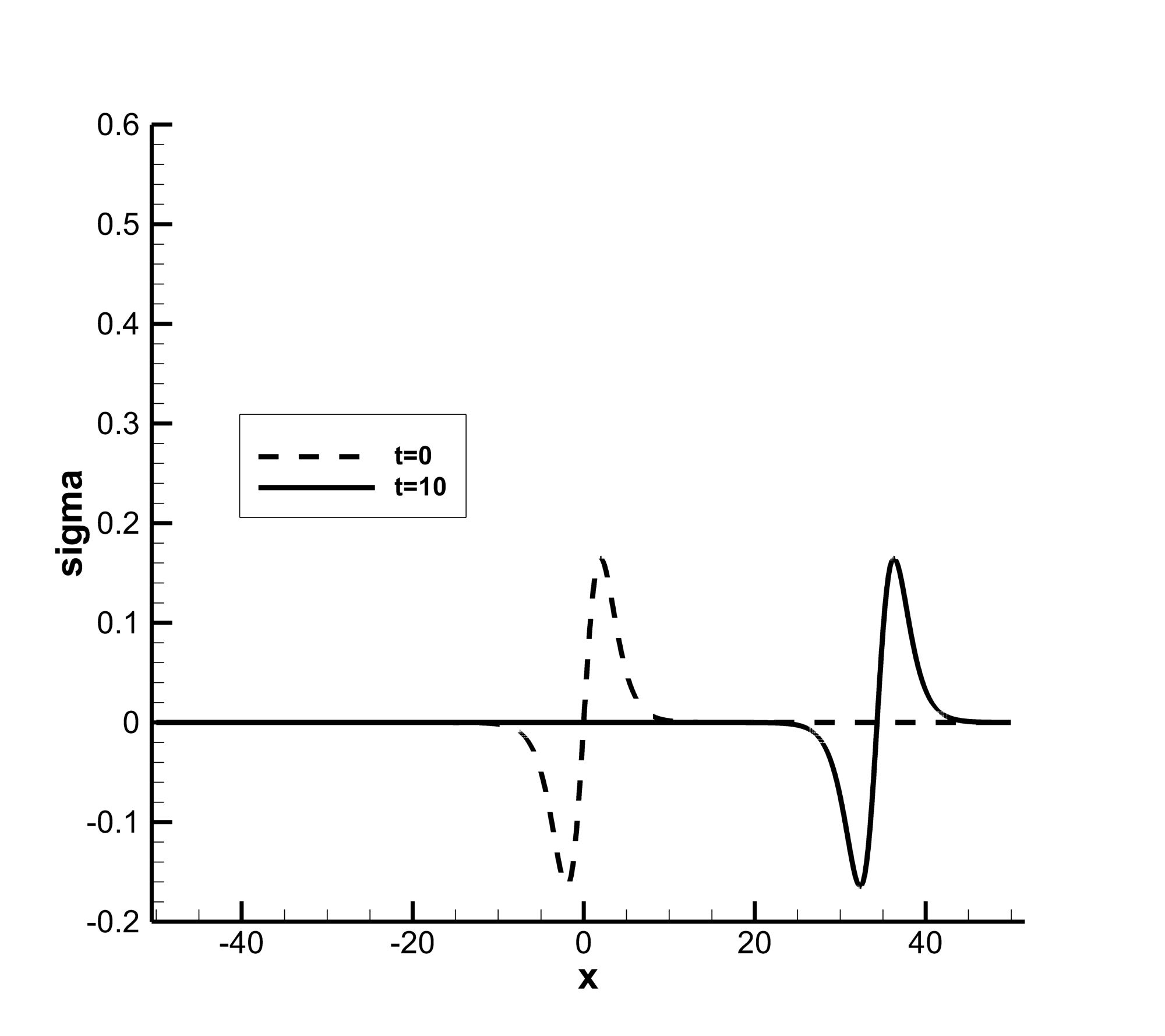}}
\end{subfigure} 

\begin{subfigure}[$\overline p$]{
     \label{fig:flat_bottom_p}
     \includegraphics[width=0.42\textwidth]{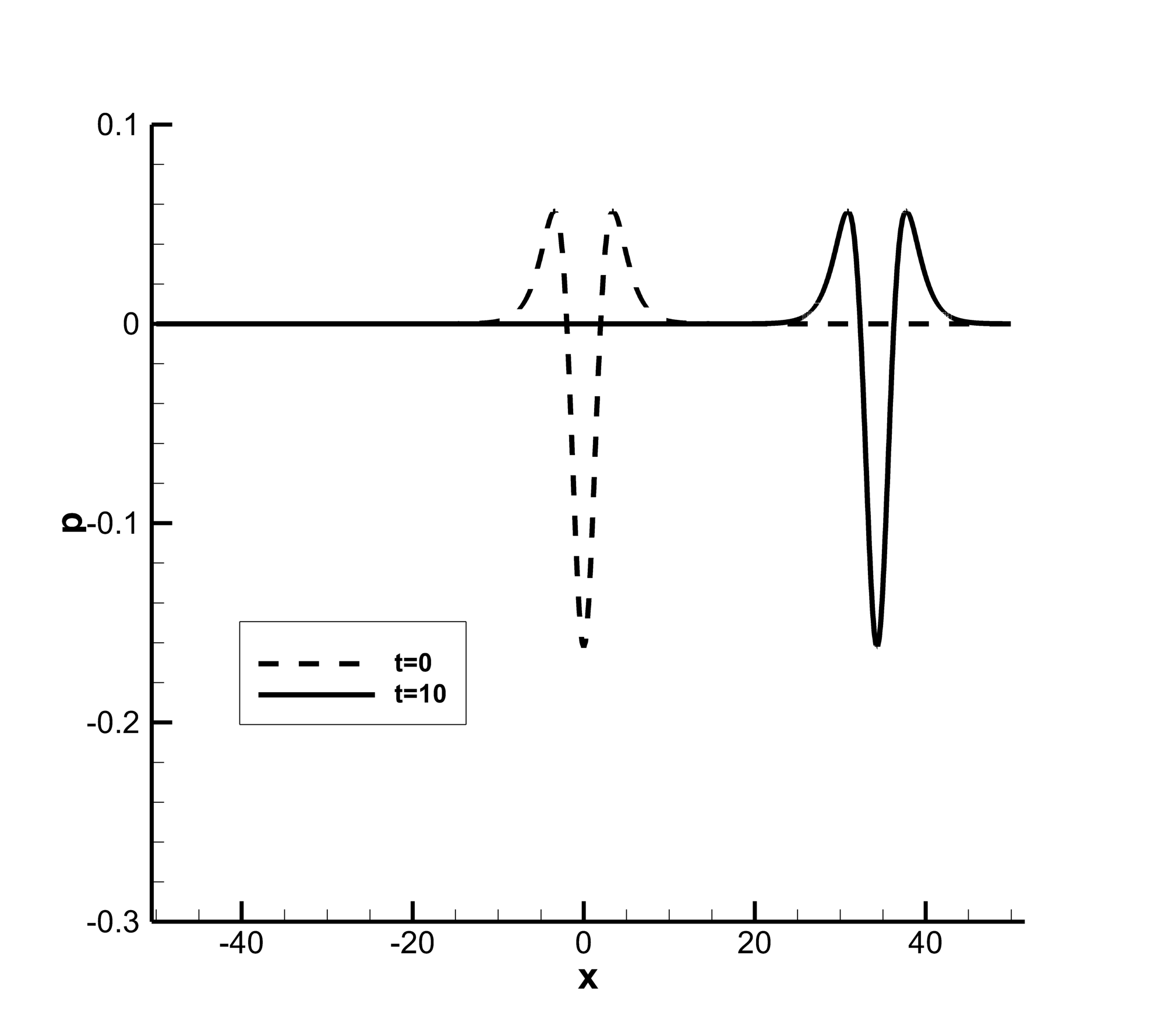}}    
\end{subfigure}
\begin{subfigure}[$\overline p_b$]{
      \label{fig:flat_bottom_pb}
      \includegraphics[width=0.42\textwidth]{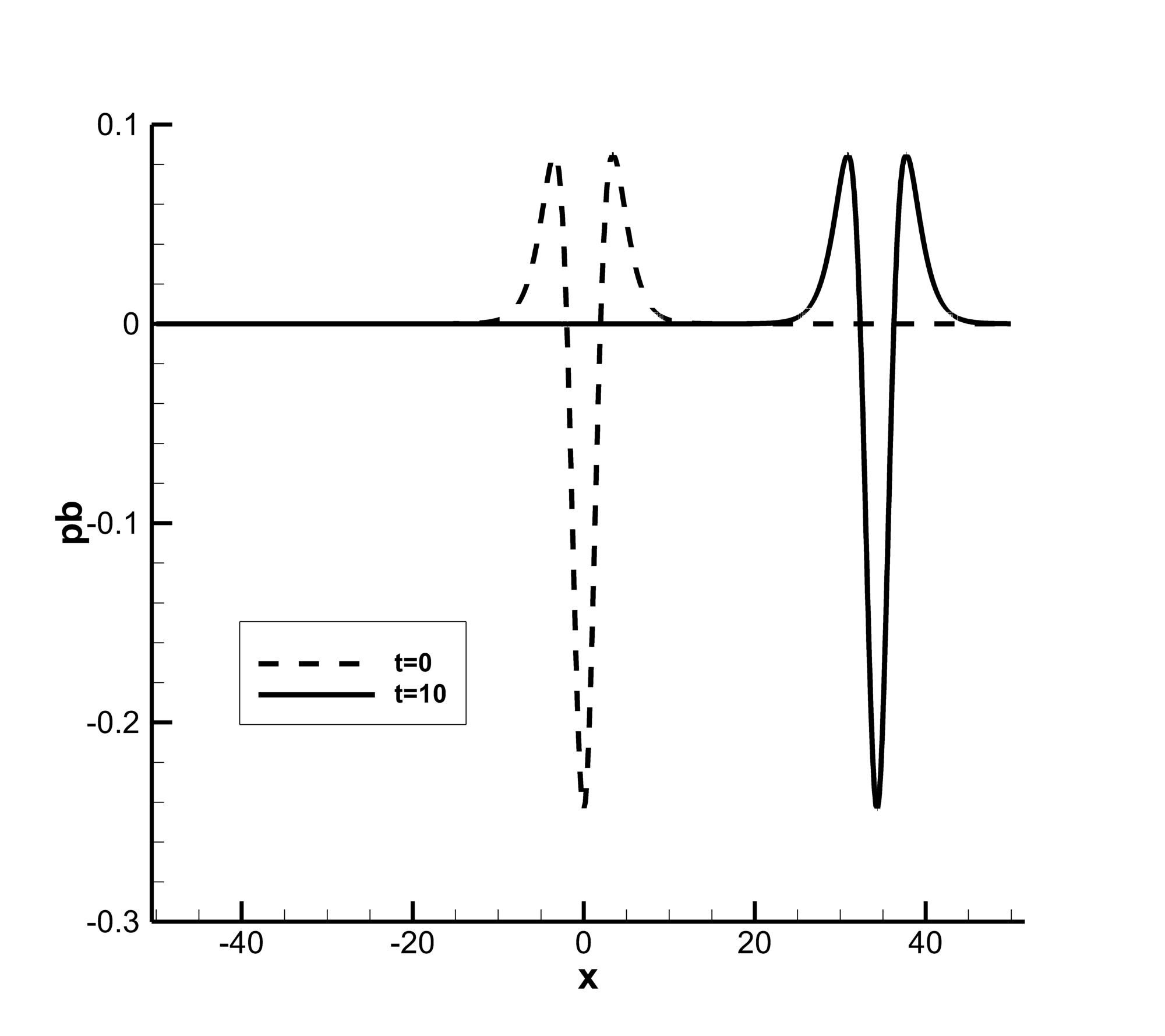}}
\end{subfigure} 
\caption{Initial condition (dashed lines) and solution at the simulation time $t_{\mathrm{end}}=10$ (solid lines) for the solitary wave over a flat bottom test case.}
\label{fig:flat_bottom}
\end{figure}

In Figure \ref{fig:flat_bottom} the initial condition (dashed lines) is shown, together with the solution at simulation time $t=10$ (solid line), obtained with a fourth order accurate scheme. We can notice that, as expected, the soliton simply propagates without changing its shape. 

\subsection{Solitary wave over a step}
The test case presented in this section is the same as the one reported  in \cite{seabrasantos:1987}, where both experimental and numerical results are provided. The domain is $\Omega = [-16,17]$, while the final simulation time is \mbox{$t_{\mathrm{end}}=10.74$}. Transmissive boundary conditions have been imposed. 
The initial amplitude of the soliton is $A=0.0365$, the still water depth is $H=0.2$, and at the  the initial time the soliton is 
centered at $x_0=-3$. The obstacle is a step of height $\Delta H_{\mathrm{obs}}=0.1$, located at $x_{\mathrm{obs}}=0$. With respect to the experiment presented in \cite{seabrasantos:1987}, the step has been smoothed out using the error function
$ z_b(x) = 0.05 \, ( \, \erf(8 x) + 1 \,). $   

\begin{figure}
\centering
\begin{subfigure}[$t=0$]{
     \label{fig:step2a}
     \includegraphics[width=0.42\textwidth,trim=10 10 50 50,clip]{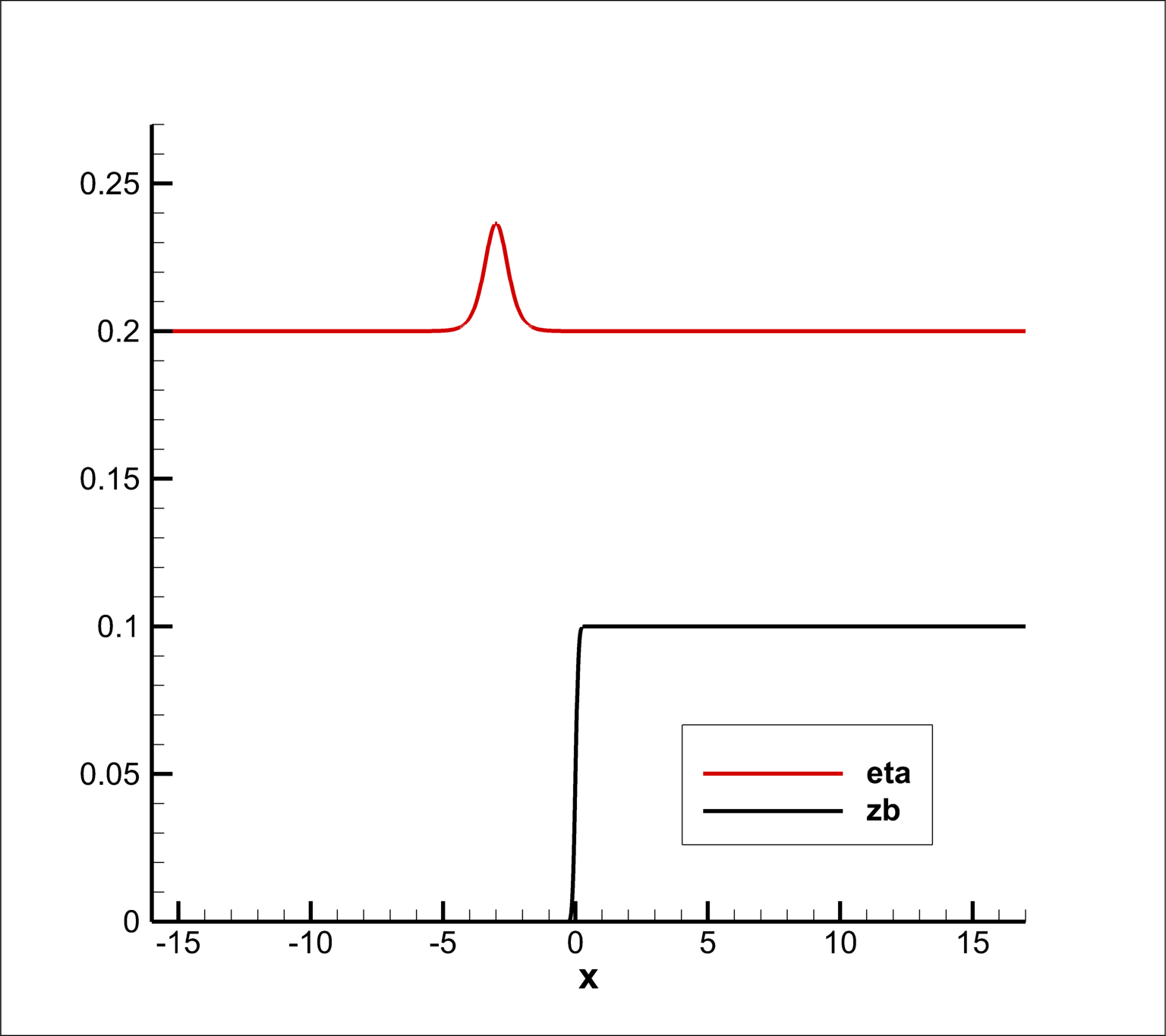}}    
\end{subfigure}
\begin{subfigure}[$t=2$]{
      \label{fig:step2b}
     \includegraphics[width=0.42\textwidth,trim=10 10 50 50,clip]{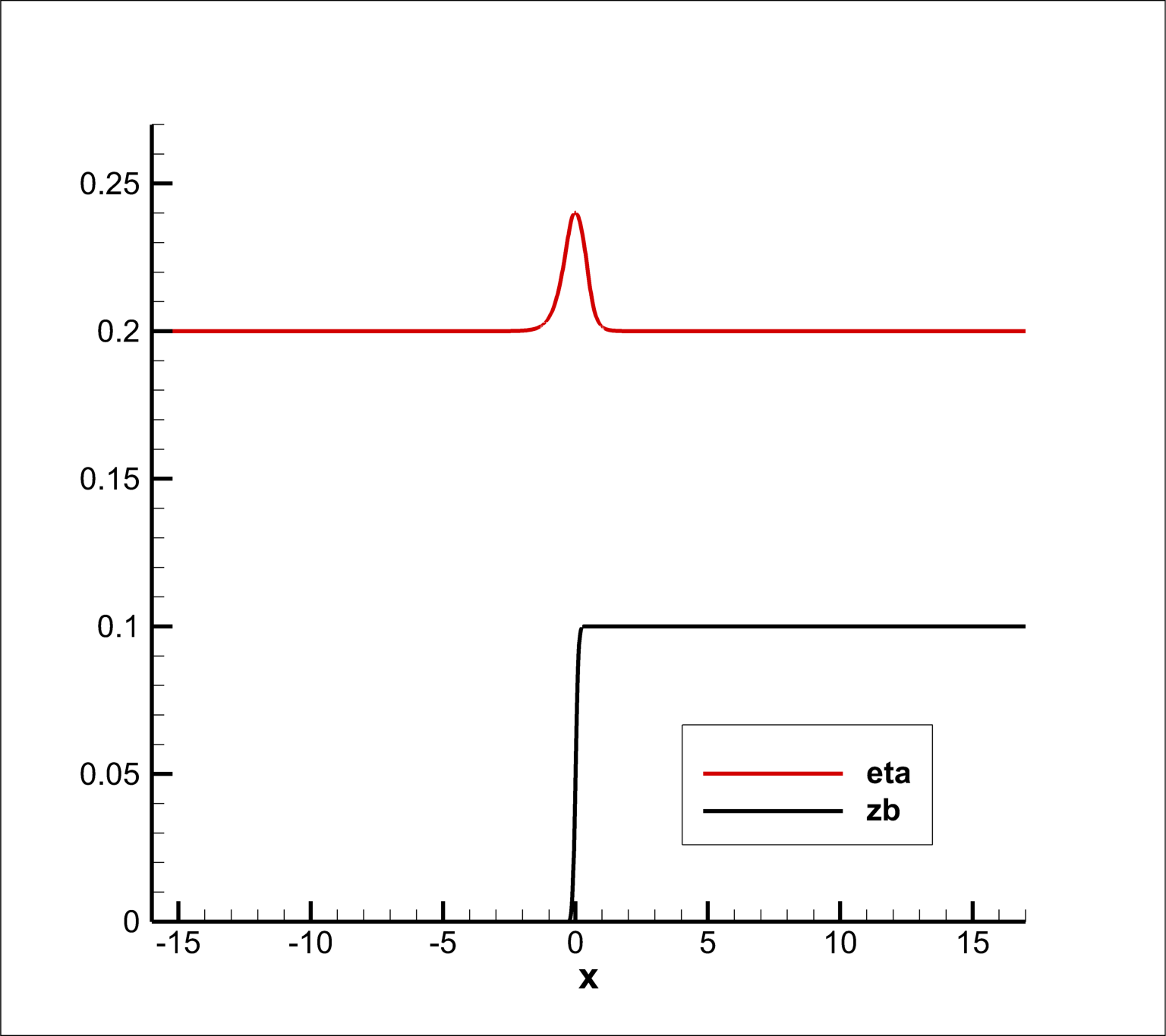}}
\end{subfigure} 

\begin{subfigure}[$t=4$]{
     \label{fig:step2c}
    \includegraphics[width=0.42\textwidth,trim=10 10 50 50,clip]{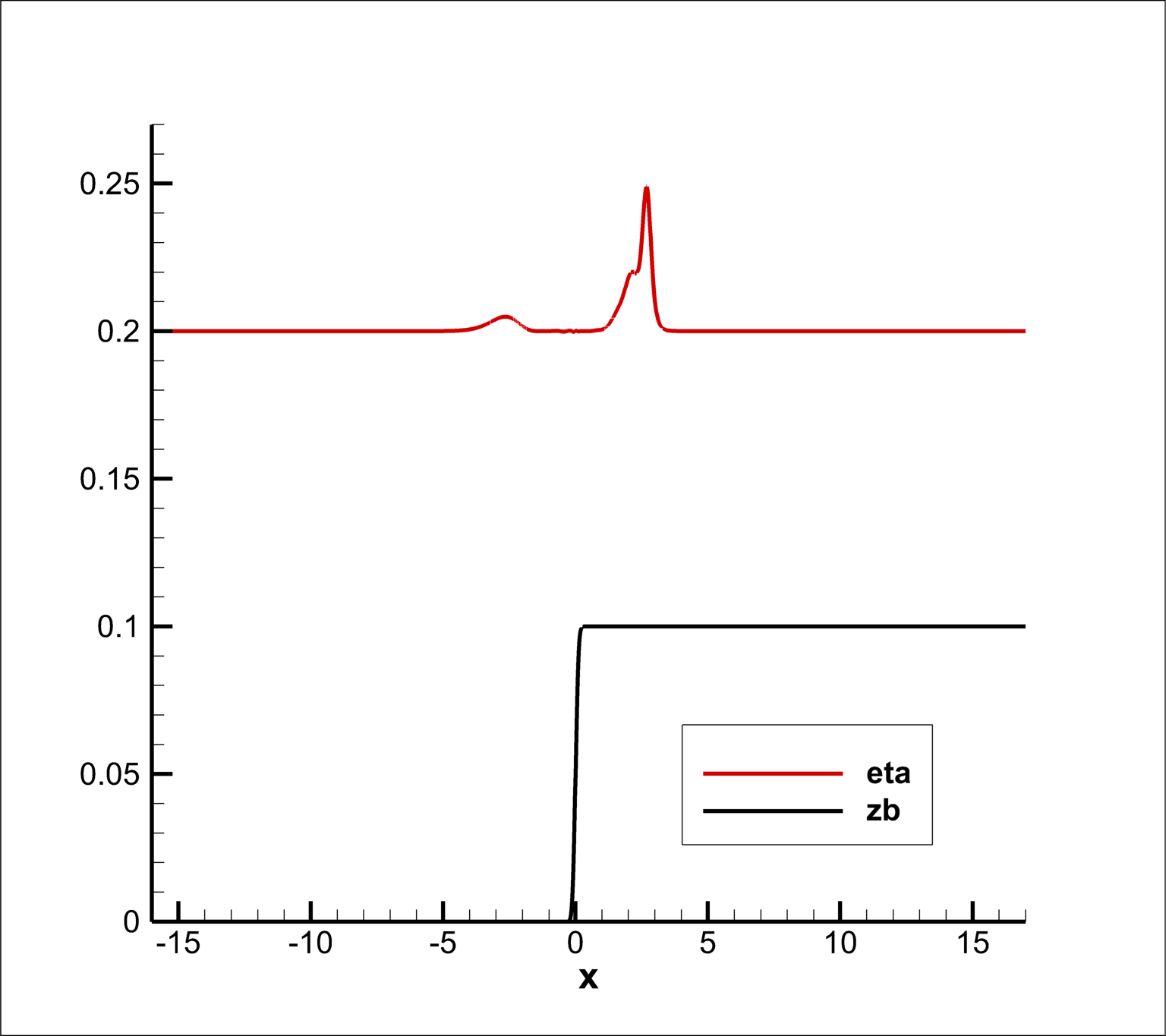}}    
\end{subfigure}
\begin{subfigure}[$t=6$]{
      \label{fig:step2d}
     \includegraphics[width=0.42\textwidth,trim=10 10 50 50,clip]{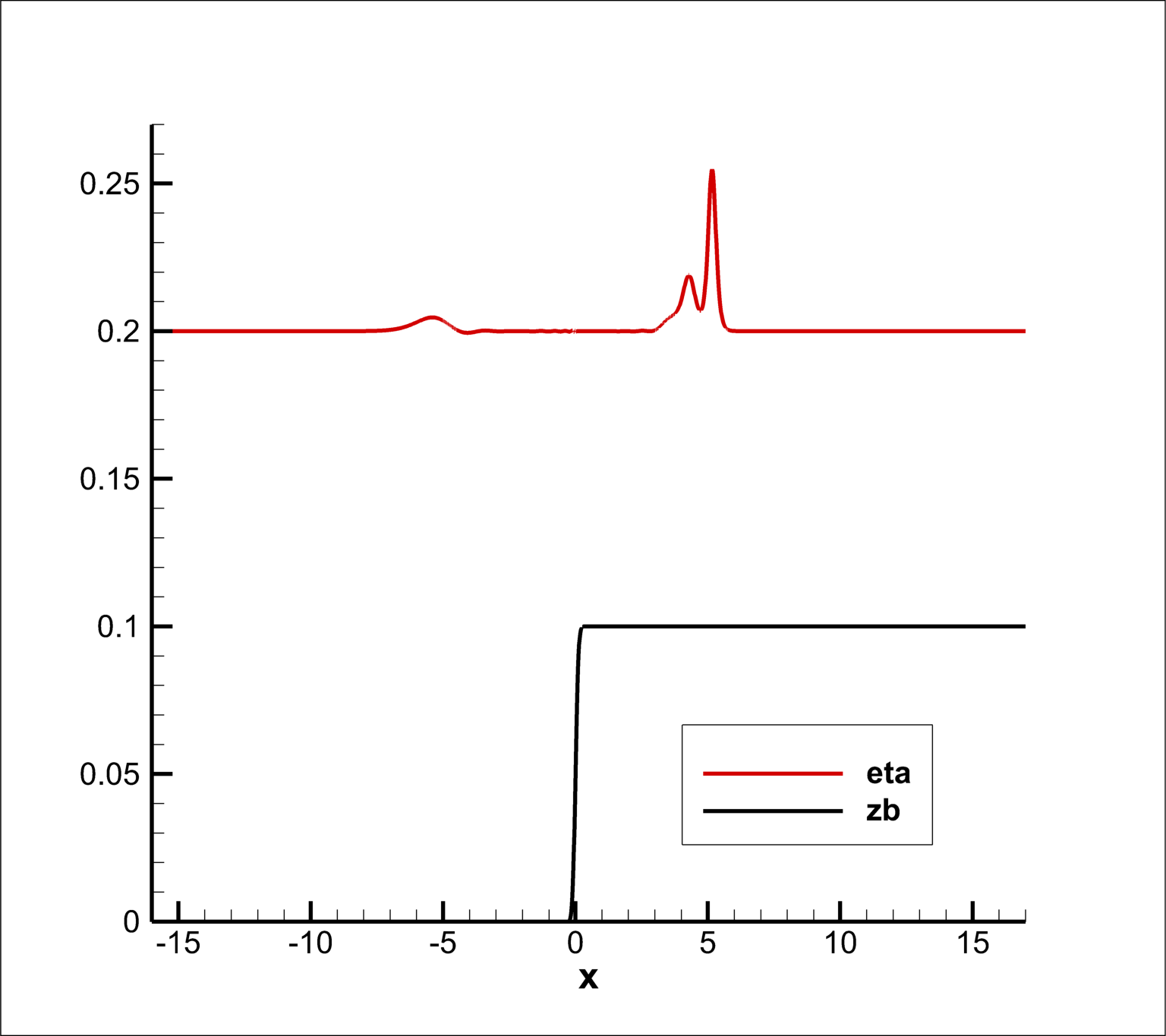}}
\end{subfigure} 

\begin{subfigure}[$t=8$]{
     \label{fig:step2e}
     \includegraphics[width=0.42\textwidth,trim=10 10 50 50,clip]{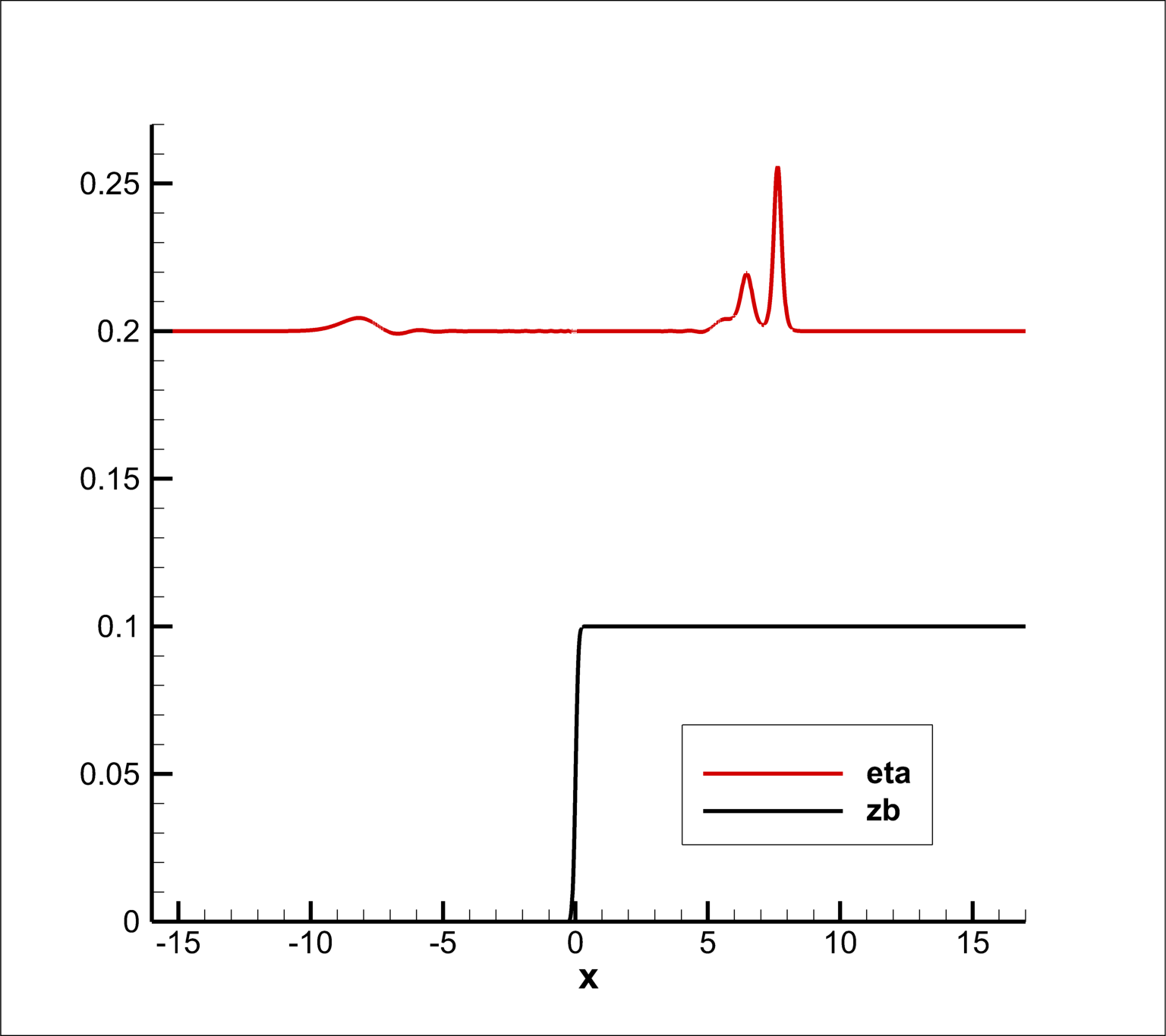}}    
\end{subfigure}
\begin{subfigure}[$t=10$]{
      \label{fig:step2f}
      \includegraphics[width=0.42\textwidth,trim=10 10 50 50,clip]{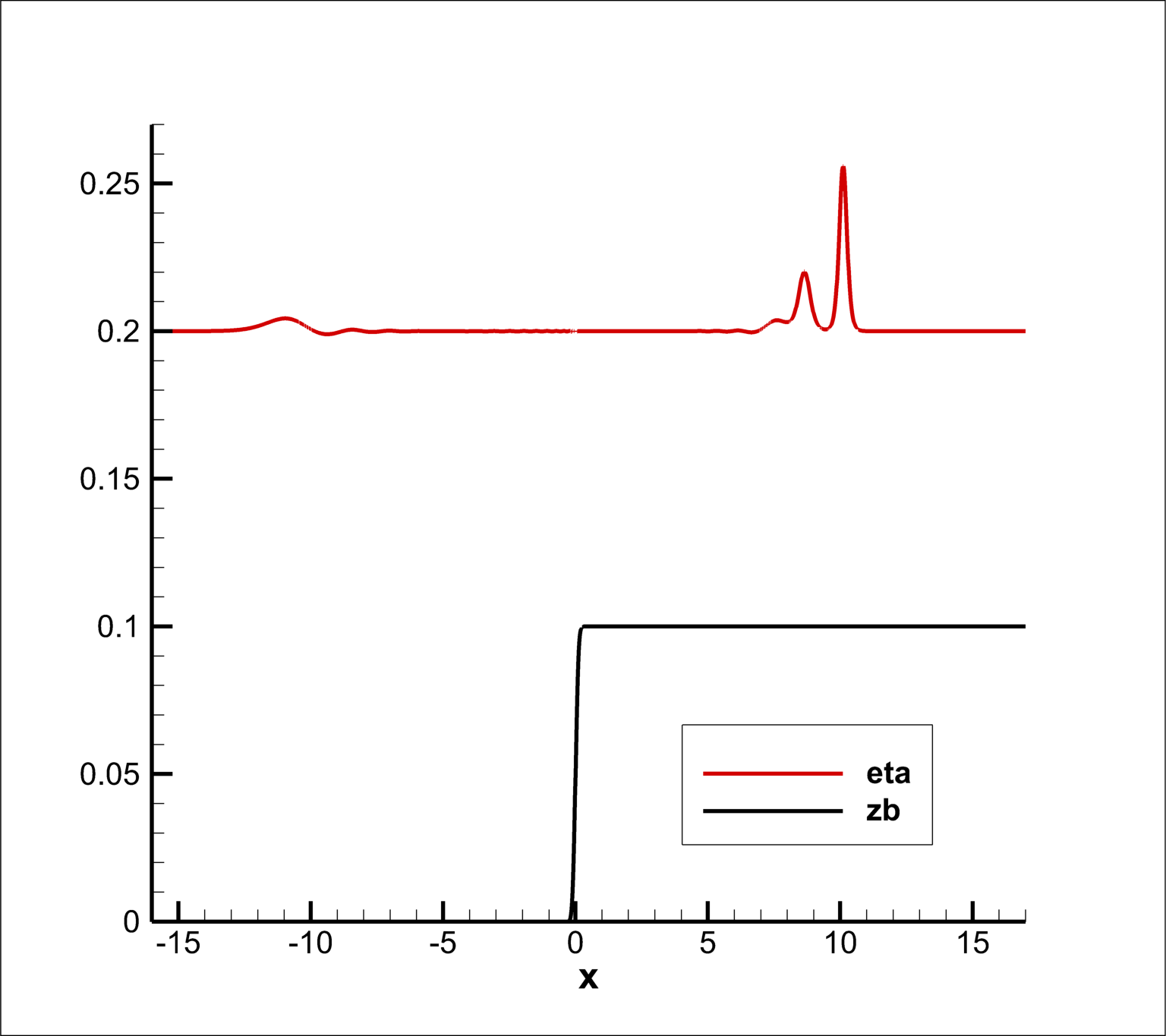}}
\end{subfigure} 
\caption{Snapshots of the free-surface $\eta=h+z_b$ (red line) at different times, for the solitary wave over a step test case. The obstacle (black line) is also represented.}
\label{fig:step2}
\end{figure}

The results have been obtained with a grid spacing of $\Delta x = 0.0076$ and a nominally fourth order accurate scheme ($N=3$).
In Figure \ref{fig:step2} we first show  some snapshots of the free surface $\eta=h+z_b$ at different times. As pointed out in \cite{seabrasantos:1987}, after the interaction with the obstacle, the incident solitary wave first grows in amplitude, see Figure \ref{fig:step2c}, and later splits into two transmitted waves, followed by a few small dispersive waves. Together with the transmitted waves, a reflected wave arises, associated with a small train of dispersive waves, see Figures \ref{fig:step2d}-\ref{fig:step2f}. 

\begin{figure}
\centering
\begin{subfigure}[$x=-9$]{
     \label{fig:step1a}
     \includegraphics[width=0.42\textwidth,trim=10 10 50 50,clip]{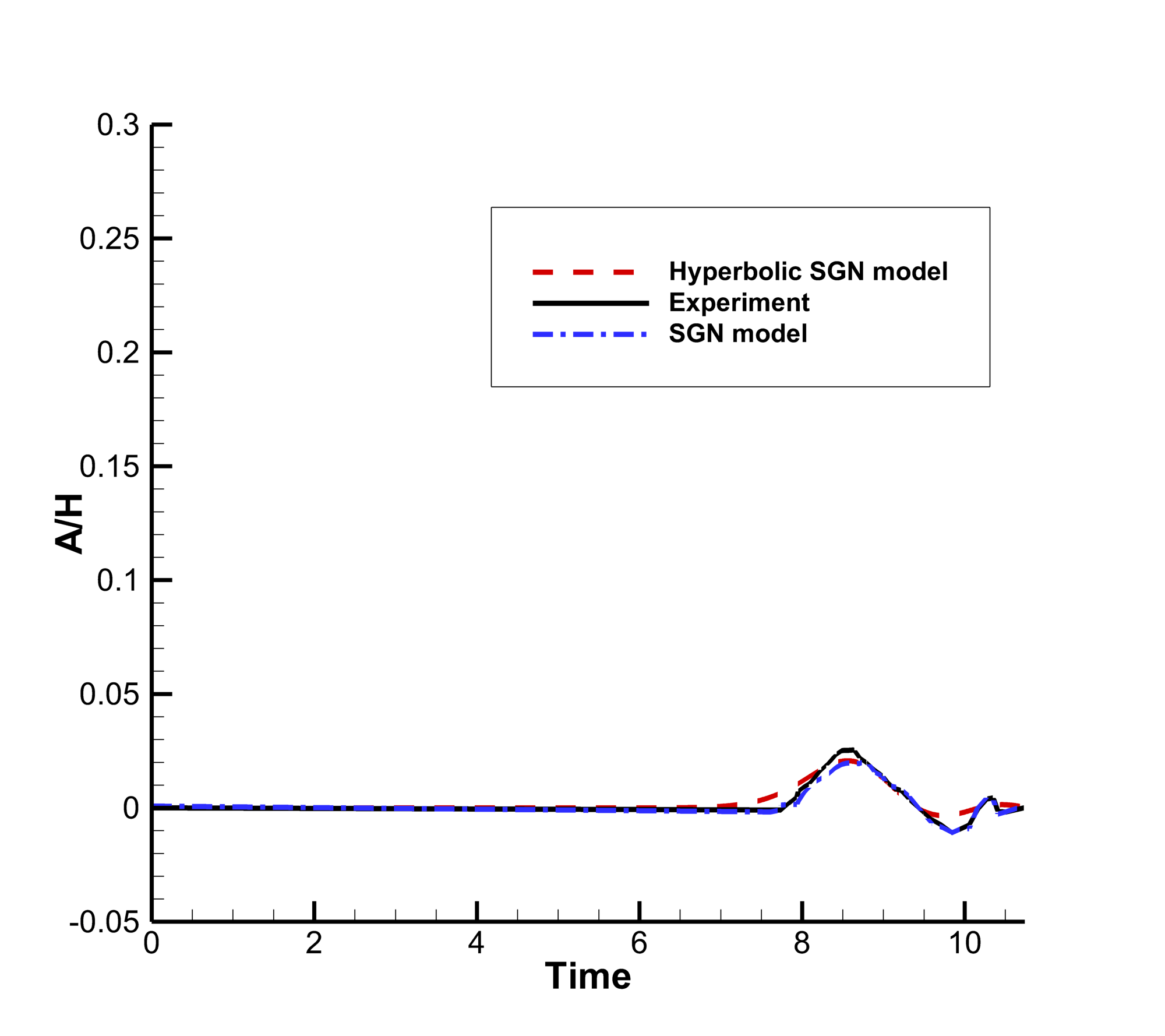}}    
\end{subfigure}
\begin{subfigure}[$x=-6$]{
      \label{fig:step1b}
     \includegraphics[width=0.42\textwidth,trim=10 10 50 50,clip]{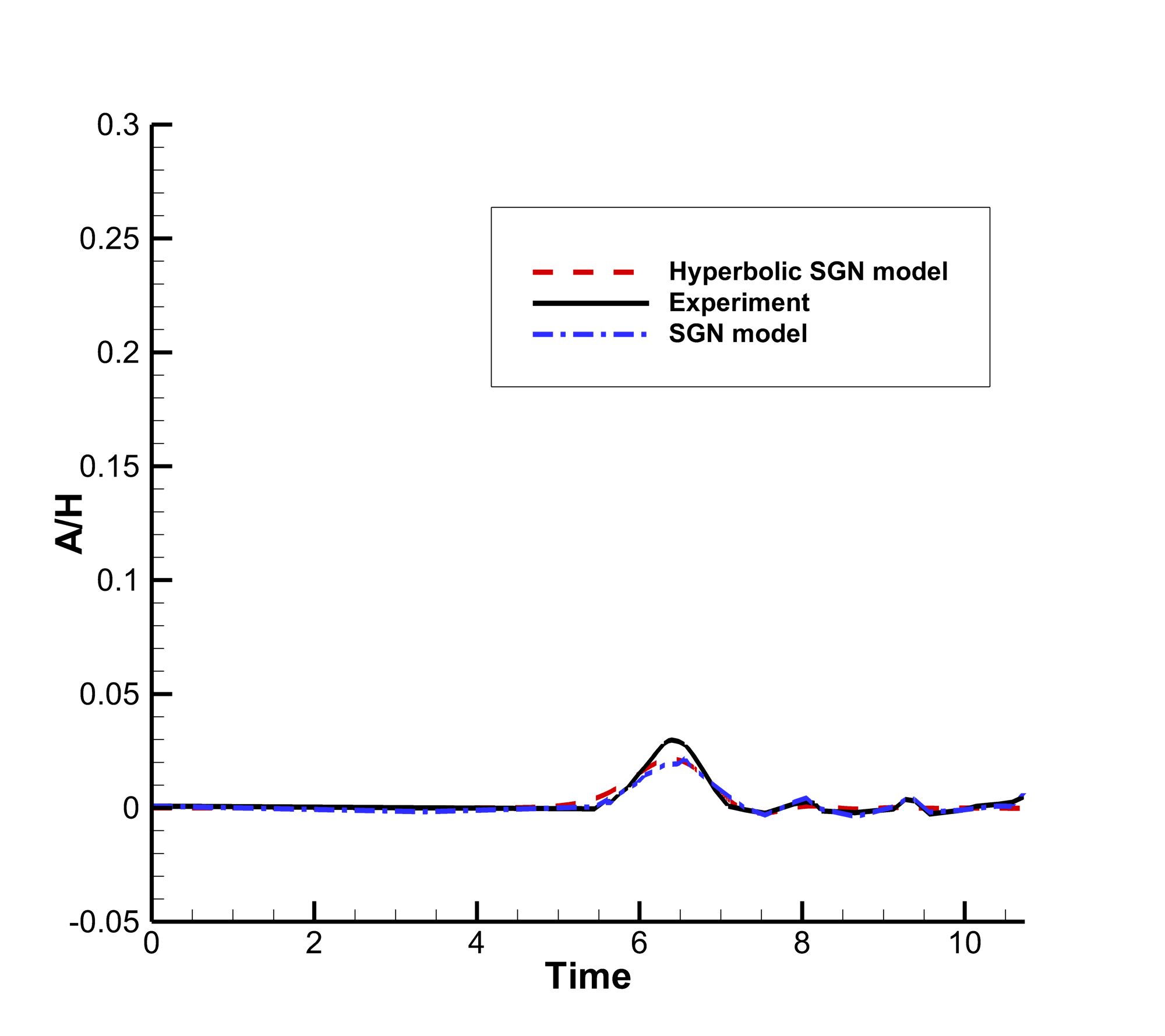}}
\end{subfigure} 

\begin{subfigure}[$x=-3$]{
     \label{fig:step1c}
    \includegraphics[width=0.42\textwidth,trim=10 10 50 50,clip]{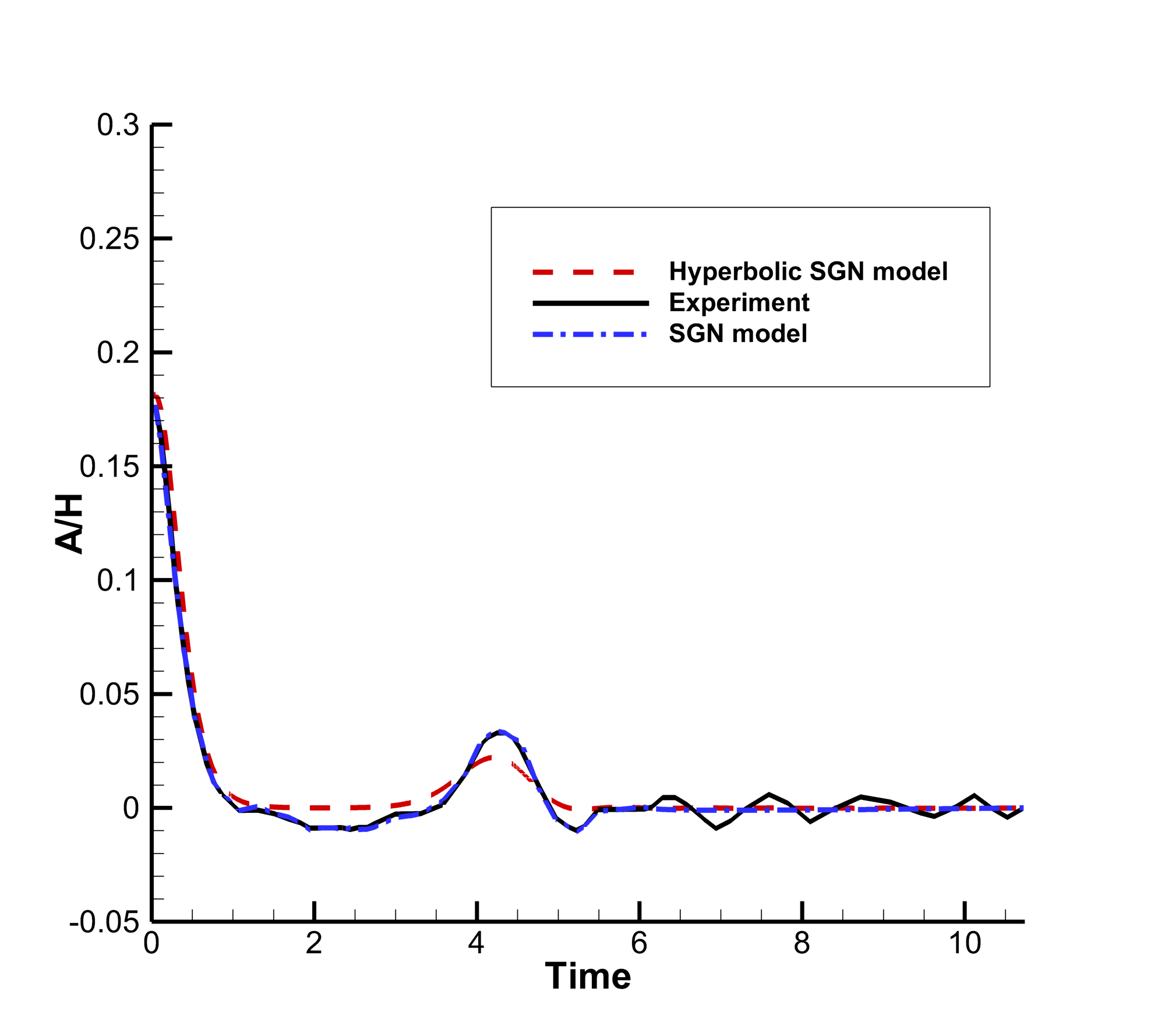}}    
\end{subfigure}
\begin{subfigure}[$x=3$]{
      \label{fig:step1d}
     \includegraphics[width=0.42\textwidth,trim=10 10 50 50,clip]{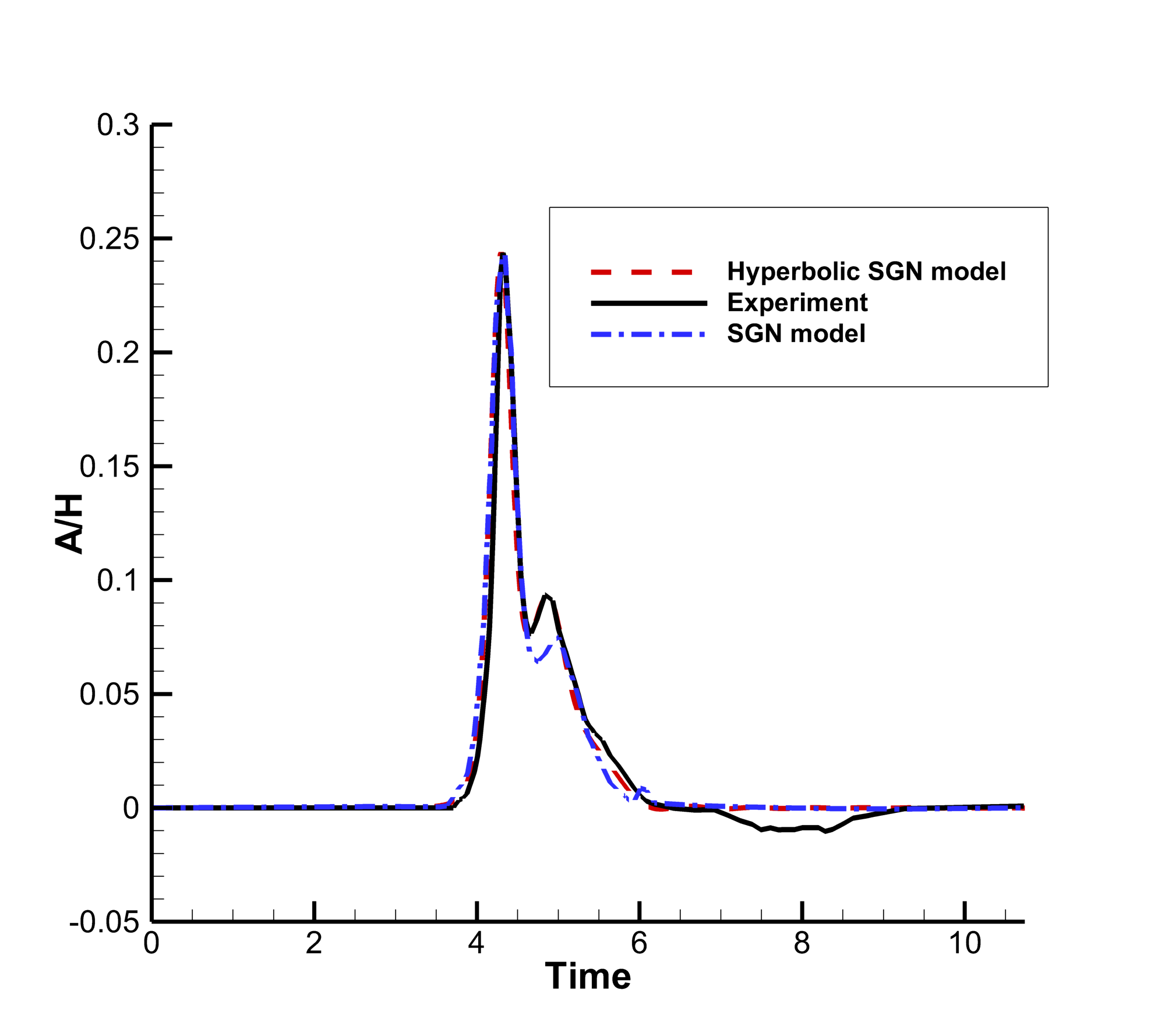}}
\end{subfigure} 

\begin{subfigure}[$x=6$]{
     \label{fig:step1e}
     \includegraphics[width=0.42\textwidth,trim=10 10 50 50,clip]{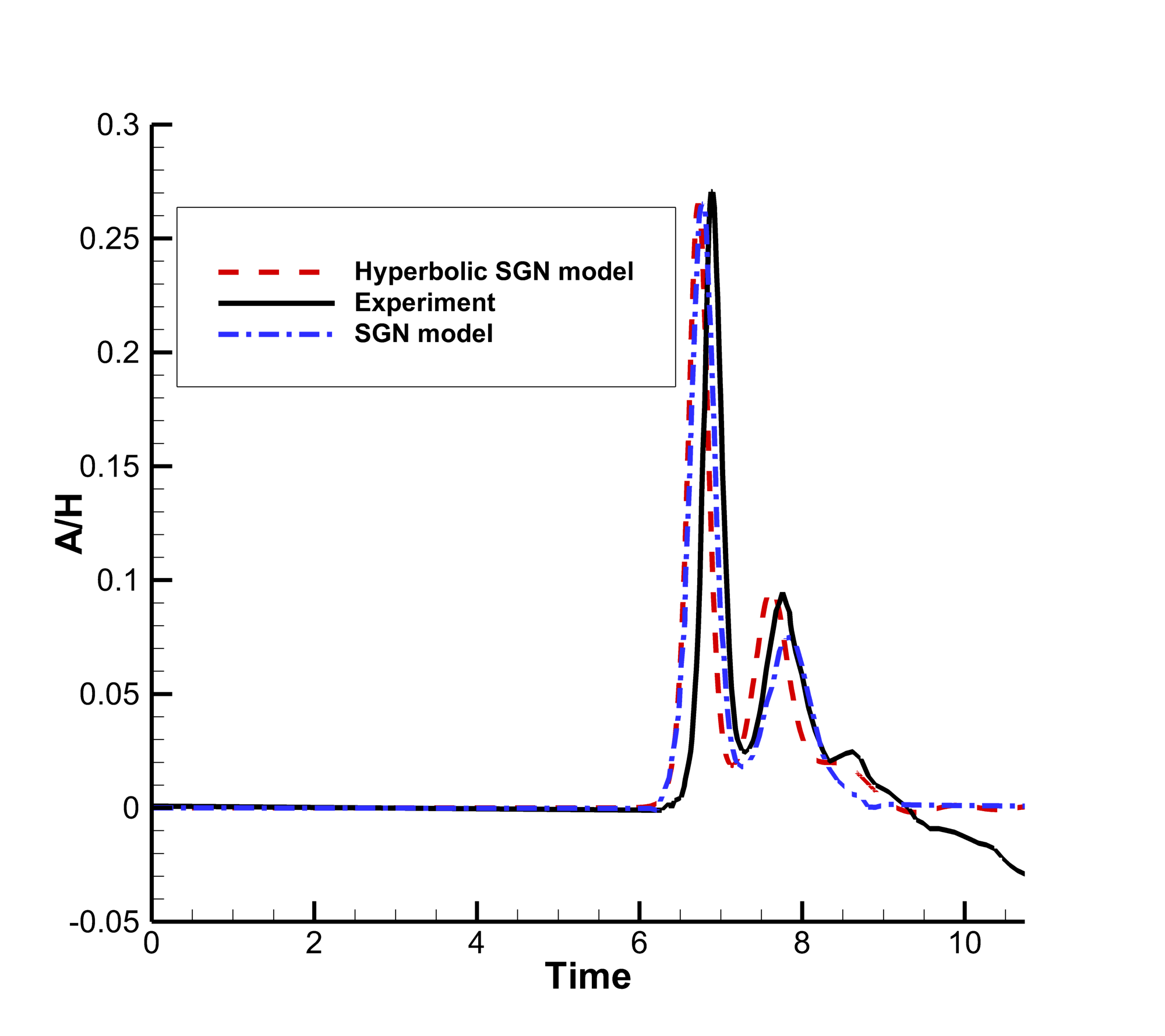}}    
\end{subfigure}
\begin{subfigure}[$x=9$]{
      \label{fig:step1f}
      \includegraphics[width=0.42\textwidth,trim=10 10 50 50,clip]{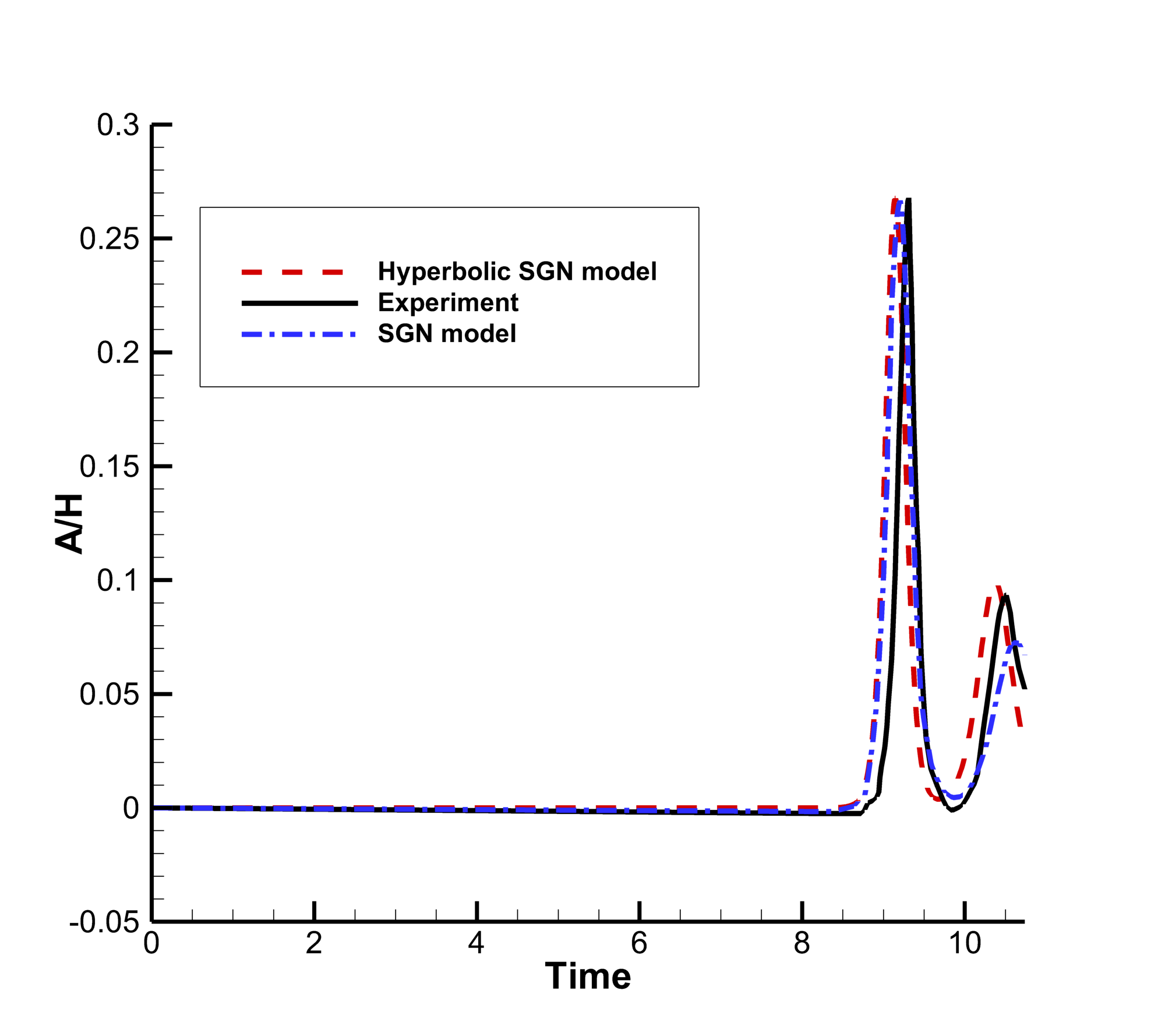}}
\end{subfigure} 
\caption{Time evolution of the quantity $A/H$ at different locations for the solitary wave over a step test case. Black line: experimental data from \cite{seabrasantos:1987}. Blue dash dot line: numerical results from \cite{seabrasantos:1987}. Red dashed line: numerical results obtained with the hyperbolic SGN model (\ref{eq:hgnb}).}
\label{fig:step1}
\end{figure}

In Figure \ref{fig:step1}, the results obtained with model \eqref{eq:hgnb} and the experimental and numerical results presented in \cite{seabrasantos:1987} are compared with each other. Six different locations are chosen in the domain and the time evolution of the variable $A/H$, which is the ratio of the wave amplitude and the still water depth, is considered. It should be noticed that, since the initial position of the soliton is $x=-3$, Figures \ref{fig:step1a}, \ref{fig:step1b} and \ref{fig:step1c} show the reflected waves, while Figures \ref{fig:step1d}, \ref{fig:step1e} and \ref{fig:step1f} show the transmitted waves. We can observe a good agreement between our results and both the experimental and numerical results of \cite{seabrasantos:1987}, with the exception of a small phase shift and some differences in the wave amplitudes. However we can notice that, in our case, the correct amplitude of the trasmitted waves is better reproduced and, moreover, a third small transmitted wave, which is completely absent in the numerical results presented in \cite{seabrasantos:1987}, is correctly captured, see Figures \ref{fig:step1d}, \ref{fig:step1e}, \ref{fig:step1f}. 

\begin{figure}
\centering
\begin{subfigure}[]{
     \includegraphics[width=0.47\textwidth,trim=10 10 50 50,clip]{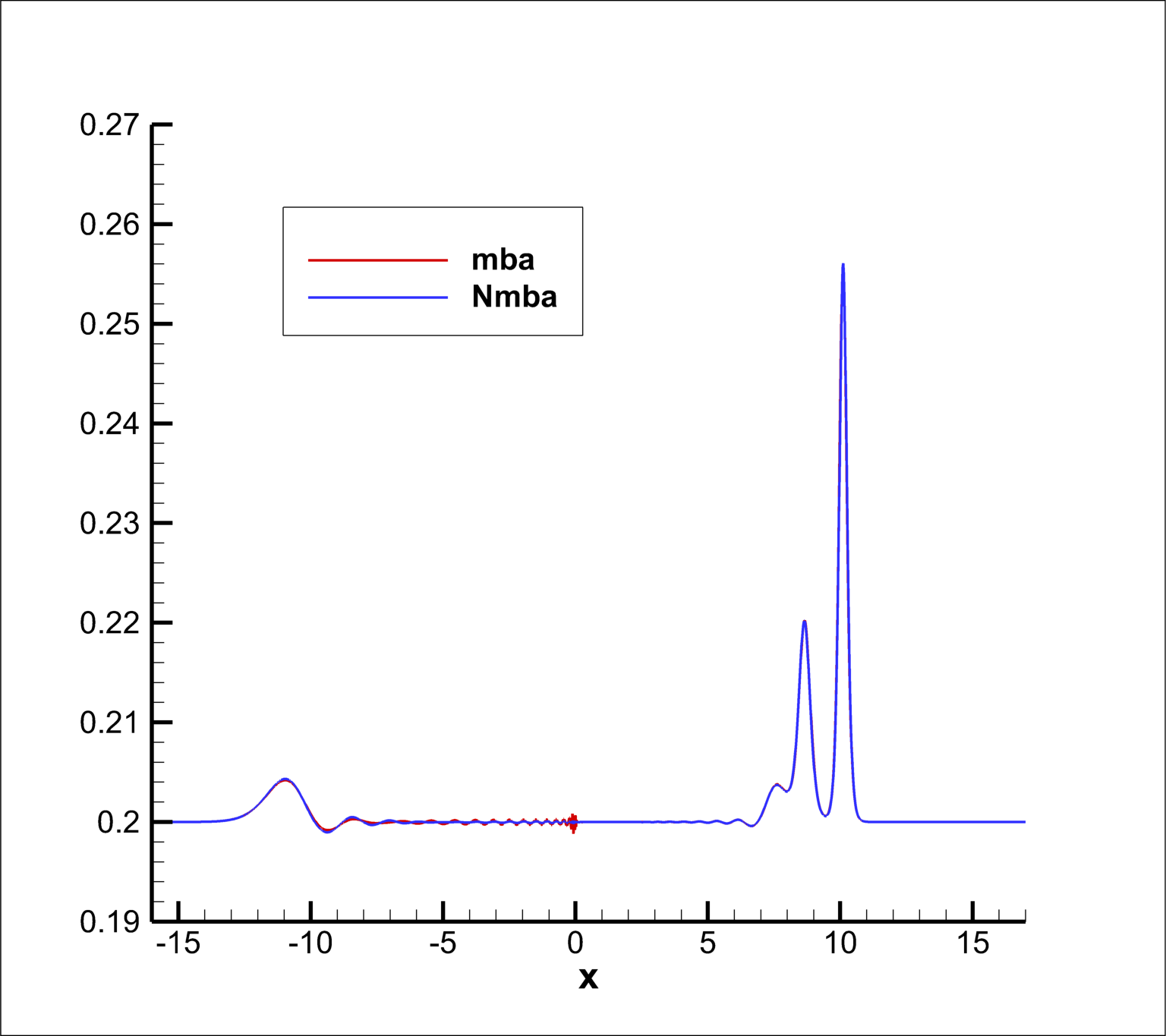}}    
\end{subfigure}
\begin{subfigure}[]{
     \includegraphics[width=0.47\textwidth,trim=10 10 50 50,clip]{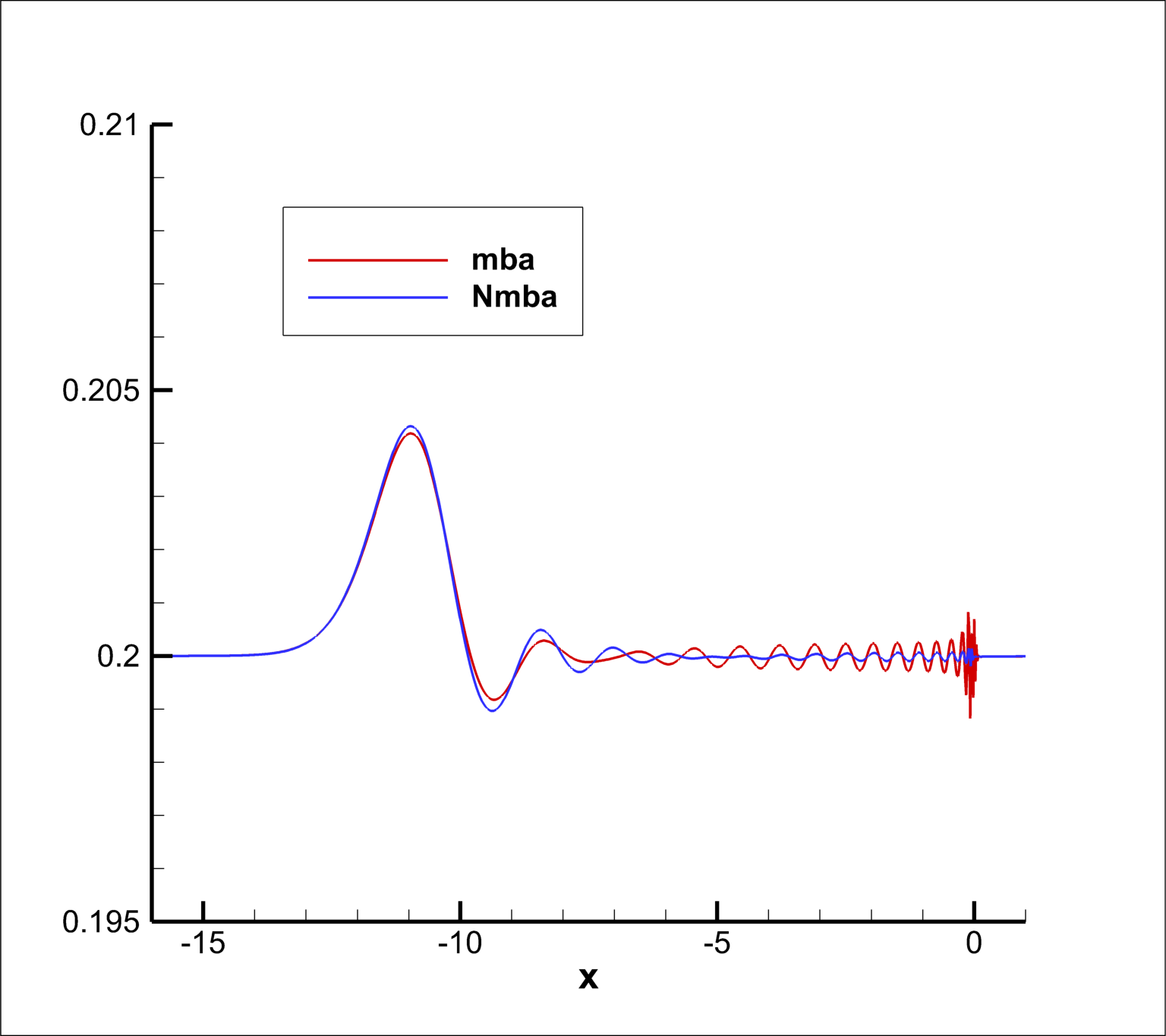}}
\end{subfigure} 
\caption{Snapshots of the free-surface $\eta=h+z_b$ at $t=10$, for the solitary wave over a step test case. Blue line: model \eqref{eq:hgnb}, i.e. hyperbolic reformulation of the SGN model without the mild bottom approximation. Red line: model \eqref{eq:hgnb_mildbot}, hyperbolic reformulation of the SGN model with the mild bottom approximation. Figure (b) is obtained by zooming in Figure (a).}
\label{fig:hsm-echgnb}
\end{figure}

Since the present test case involves a strongly varying topography, we have realized an additional simulation employing the simplified model (\ref{eq:hgnb_mildbot}) with mild bottom approximation, in order to investigate the differences with respect to the solutions of the new system \eqref{eq:hgnb}, which does \textit{not} make the mild bottom assumption. 
In Figure \ref{fig:hsm-echgnb}, the free surface elevation is represented for both, the new model (blue line) and model \eqref{eq:hgnb_mildbot} (red line). We can see that the simpler model \eqref{eq:hgnb_mildbot}  introduces much larger spurious oscillations in correspondence to the bottom step. Notice that the same behaviour can be observed for different time instants.   
In order to verify if the oscillating behaviour of model \eqref{eq:hgnb_mildbot} is linked to the particular choice of the numerical method, we have realized an additional simulation with a classical second order finite volume scheme on a very fine mesh for both, the new hyperbolic model without mild bottom assumption proposed in this paper \eqref{eq:hgnb} and model \eqref{eq:hgnb_mildbot} proposed in \cite{escalante:2018}. Also in this case, model \eqref{eq:hgnb_mildbot} presents far more spurious oscillations in correspondence to the obstacle than the proposed new model.

\subsection{Periodic waves over a submerged bar}
\begin{figure}
\centering
\includegraphics[width=0.9\textwidth]{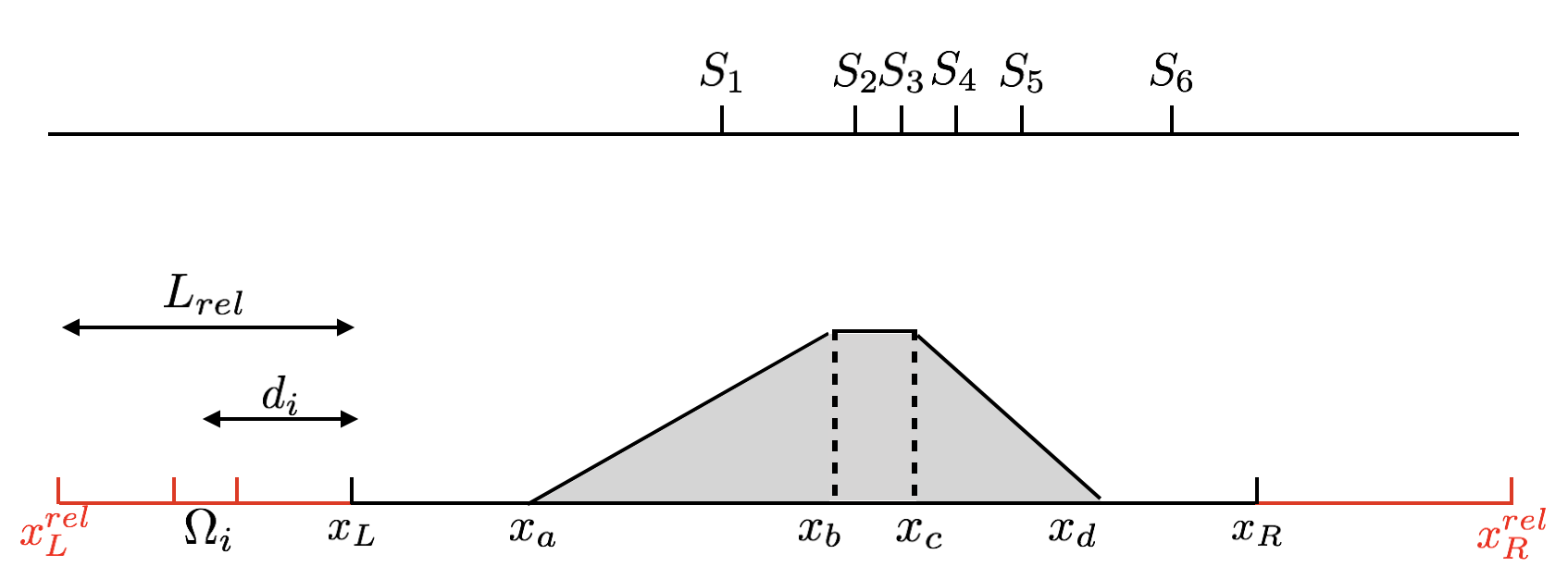}
\caption{Computational domain and wave gauges locations for the periodic waves over a submerged bar test case. Computational domain: $\Omega=[0,40]$. Relaxation zones: $[x_L^{rel}=-10,x_L=0]$, $[x_R=40,x_R^{rel}=50]$. Wave gauges locations: $S_1=10.8$, $S_2=12.8$, $S_3=13.8$, $S_4=14.8$, $S_5=16$, $S_6=17.6$. Obstacle dimensions: $x_a=6$, $x_b=12$, $x_c=14$, $x_d=17$. Length of the relaxation zone: $L_{rel}=10$.}
\label{fig:domain_subbar}
\end{figure}

In the present section, the numerical results obtained with the new model \eqref{eq:hgnb} are compared with the experimental data reported in \cite{beji:1994}. 
A sketch of the computational domain and the bottom topography $z_b$ is presented in Figure \ref{fig:domain_subbar}, together with the position of the wave gauges ($S_i$, $i=1,\cdots,6$) at which the time evolution of the solution is provided by the experimental data. A mesh of $N_x=1200$ elements and polynomial degree $N=3$ are employed. 
Concerning boundary conditions, experimental data at an additional wave gauge $S_0$, positioned before the beginning of the obstacle, are used to impose the pattern at the inflow boundary. The experimental data provides the time evolution of the wave amplitude $A^{*}(t)$, from which the time evolution of the water depth can be readily obtained as $h^*=A^*(t)+H$, where $H=0.4$. Analogously to \cite{madsen:2003, escalante:2018}, the remaining variables have the following expressions:

\begin{equation}
u^*(t) = \frac{\sqrt{gH}A^*(t)}{h^*(t)}, \quad w^*(t) = \sigma^*(t) = p^*(t) = p_b^*(t) = 0. 
\end{equation}
In order to have a smooth transition between the target solution $\mathbf{u}^{*}$  and the solution $\mathbf{u}_h$ inside the computational domain, a wavemaker boundary condition is applied. To this aim, a relaxation length $L_{rel}$ and a corresponding relaxation zone, located outside with respect to the computational domain, are introduced (see Figure \ref{fig:domain_subbar}). If an element $\Omega_i$ falls into the relaxation zone, we compute the auxiliary parameter

\begin{equation}
m_i = \sqrt{1-\left(\frac{d_i}{L_{rel}}\right)^2},
\end{equation}
where $d_i$ the distance between the barycentre of the element $\Omega_i$ and the closest boundary (see Figure \ref{fig:domain_subbar}). The solution inside the element is then redefined as

\begin{equation}
\widetilde{\mathbf u}_h = m_i \mathbf u_h + (1-m_i)\mathbf u^*.
\end{equation}
Notice that the same procedure is employed to impose an absorbing boundary condition at the right boundary, in order to prevent wave reflection. The target solution for the absorbing boundary condition is

\begin{equation}
h^*(t)=H, \quad u^*(t)= w^*(t) = \sigma^*(t) = p^*(t) = p_b^*(t) = 0. 
\end{equation}

Since we do not know exactly at which interval of time the experimental data is provided and since the exact location of the obstacle inside of the domain is unknown, the numerical results for wave gauge $S_1$ are shifted in time in order to match the experimental data at the same location. Notice that the same shifting in time is maintained also for the other wave gauges, in order to make a fair comparison between experimental data and numerical results.

\begin{figure}
\centering
\includegraphics[width=1\textwidth]{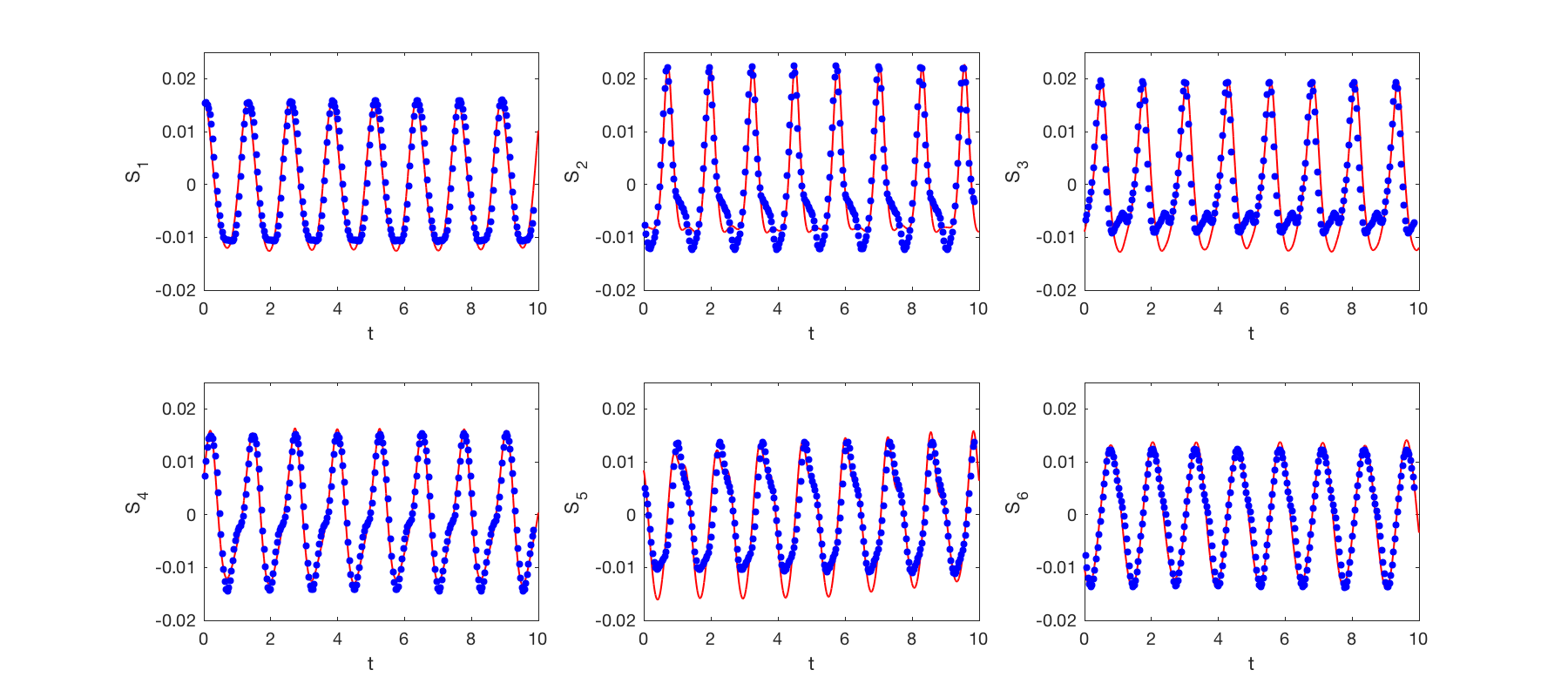}
\caption{Time evolution of the wave amplitude $A(t)$ for the periodic waves over a submerged bar test case at different wave gauges. Blue dots correspond to experimental data while red lines are the obtained numerical results.}
\label{fig:result_subbar}
\end{figure}

In Figure \ref{fig:result_subbar}, the numerical results for the time evolution of the free surface elevation $\eta (t)$ at different wave gauges are compared with the experimental data. In general, we notice a good agreement between numerical and experimental results concerning the wave period, even if the shape of the troughs at wave gauges $S_2$, $S_3$ and $S_5$ is not perfectly captured. We believe however that such discrepancies are normal when dealing with experimental data.

\subsection{Solitary wave over a Gaussian obstacle}
In order to show that the proposed approach works also in the two-dimensional case, we present the results of a test case in which a two-dimensional solitary wave impinges on a Gaussian obstacle. The computational domain is  $\Omega=[-5,35]\times[-20,20]$, to build the mesh, $200$ elements have been employed in both $x$ and $y$ directions and the polynomial degree was set to $N=3$. The Gaussian obstacle is centred at $\left( x_{\mathrm{obs}},y_{\mathrm{obs}}\right) =\left( 0,0\right) $ and its shape is

\begin{equation}
z_b(x,y) = A_g \exp{\left(-\frac{x^2+y^2}{2 \sigma_g^2}\right)},
\end{equation}
where $A_g=0.1$ and $\sigma_g = 1$. 
Concerning the soliton parameters, its initial position is $x_0=-3$, while its amplitude is $A_i = 0.1$. The still water depth is $H=0.25$ and a final integration time of $t_{\mathrm{end}}=12$ is reached. 

\begin{figure}
\centering
\begin{subfigure}[$t=0$]{
     \label{fig:gausssolitona}
     \includegraphics[width=0.47\textwidth,trim=10 10 10 50,clip]{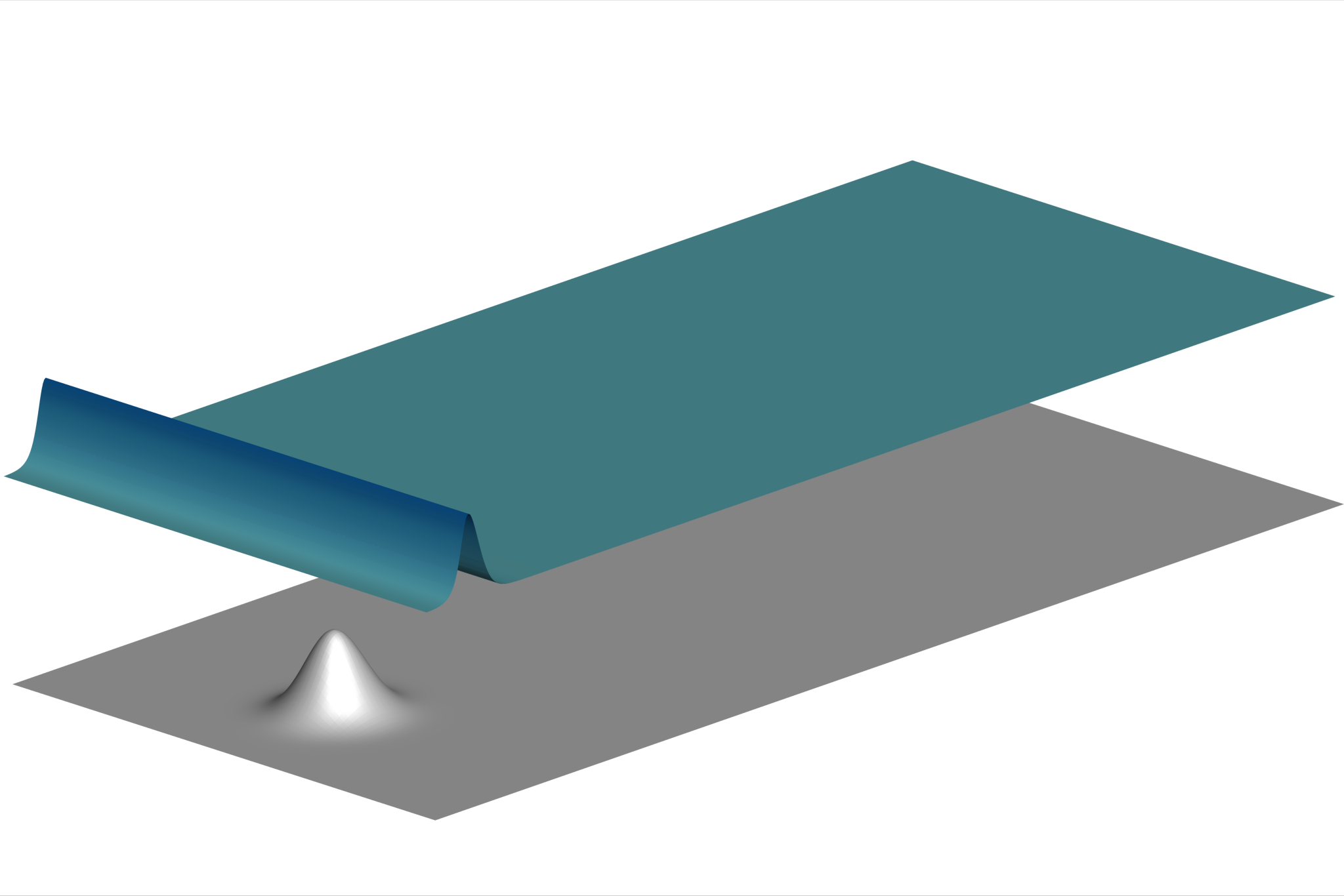}}    
\end{subfigure}
\begin{subfigure}[$t=2$]{
      \label{fig:gausssolitonb}
     \includegraphics[width=0.47\textwidth,trim=10 10 10 50,clip]{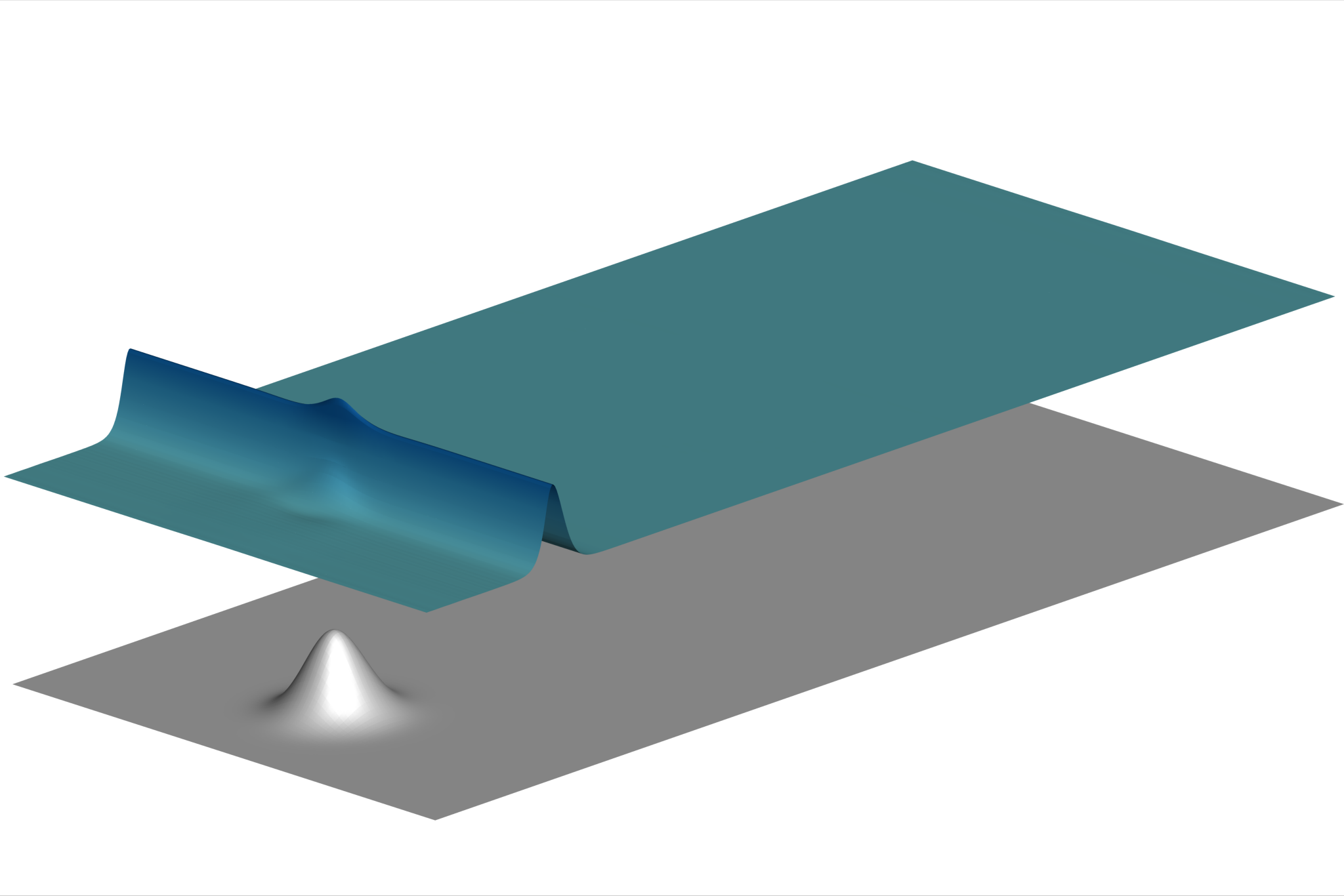}}
\end{subfigure} 

\begin{subfigure}[$t=5$]{
     \label{fig:gausssolitonc}
    \includegraphics[width=0.47\textwidth,trim=10 10 10 50,clip]{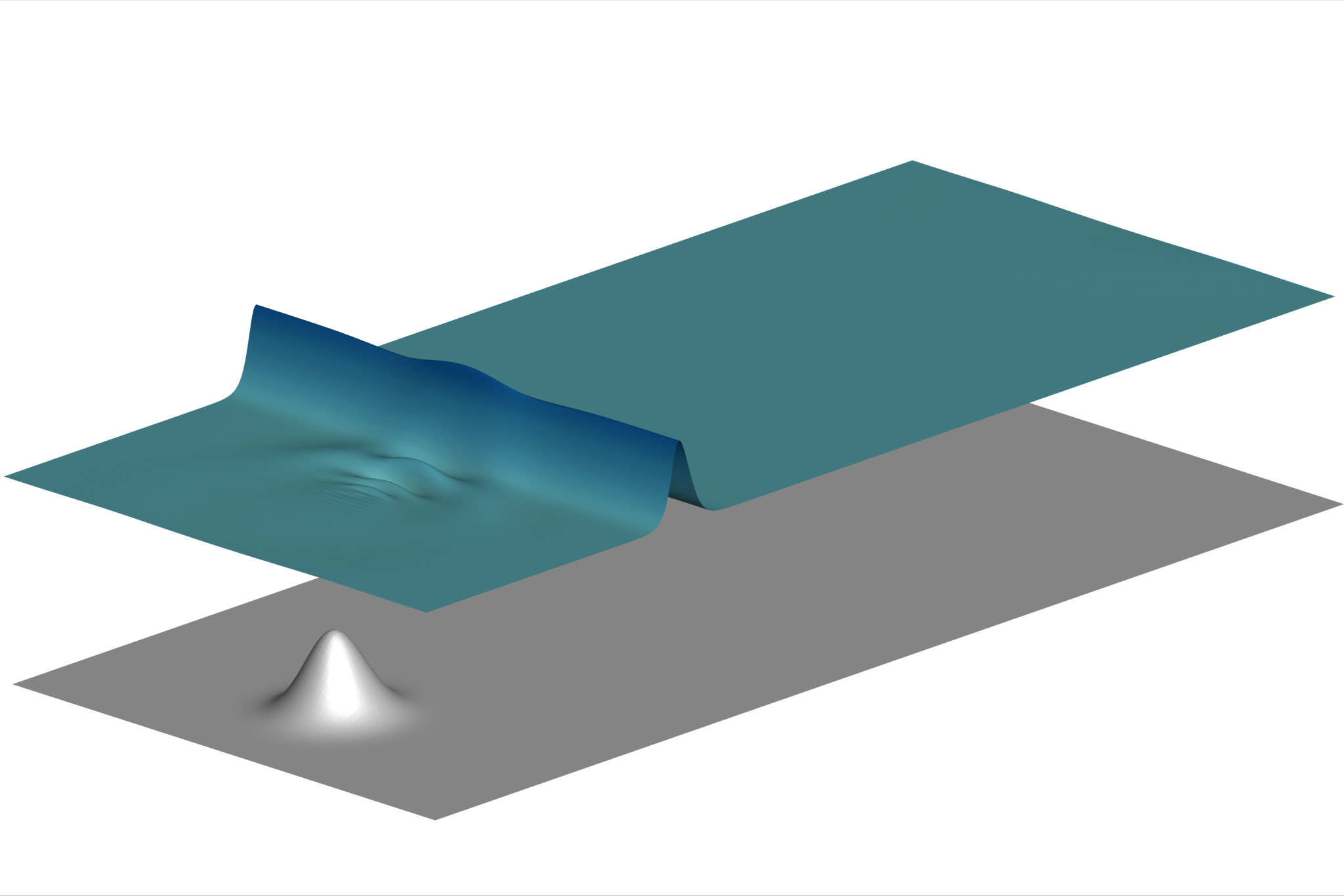}}    
\end{subfigure}
\begin{subfigure}[$t=12$]{
      \label{fig:gausssolitond}
     \includegraphics[width=0.47\textwidth,trim=10 10 10 50,clip]{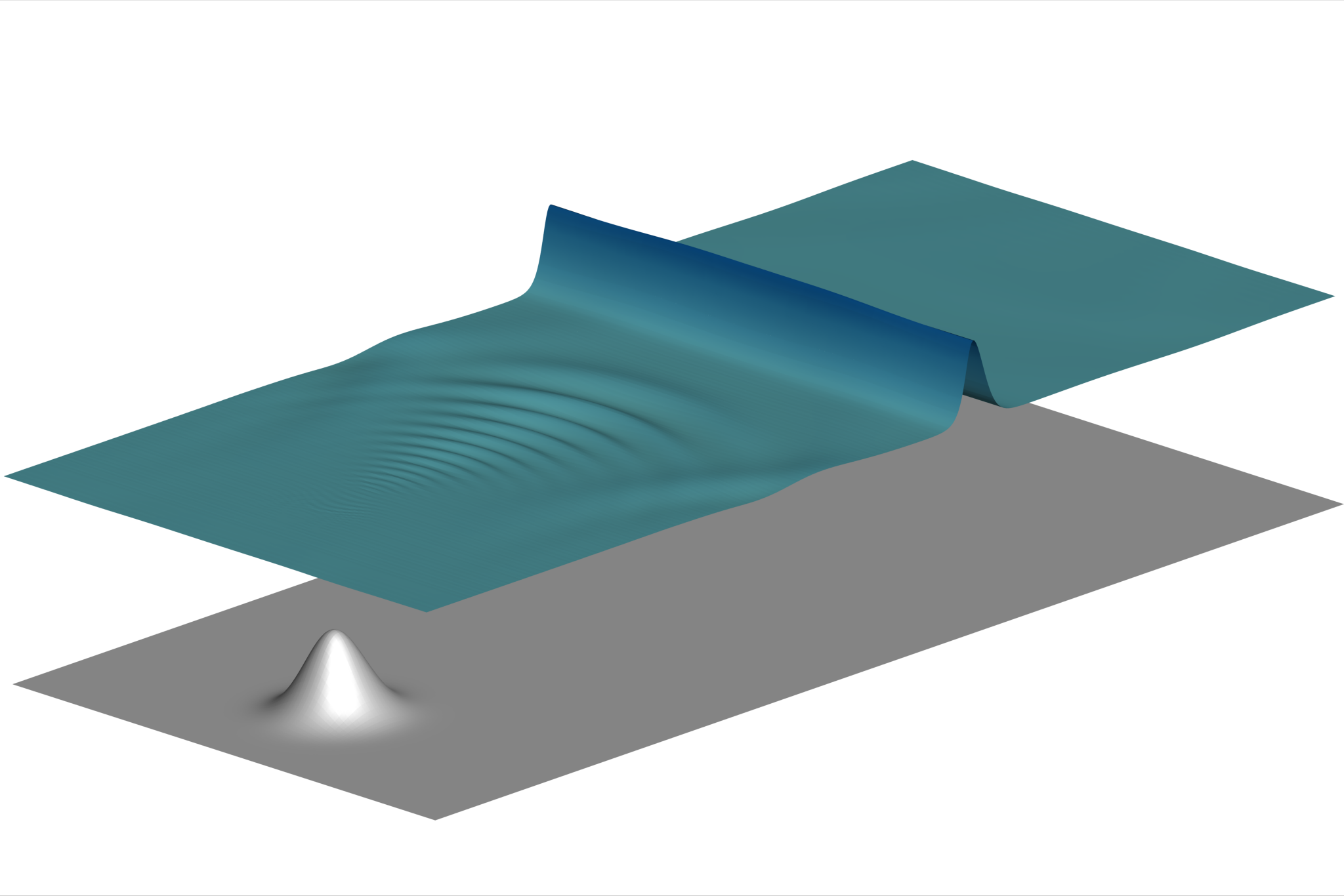}}
\end{subfigure} 
\caption{Snapshots of the free-surface $\eta=h+z_b$ (in blue) at different times, for the solitary wave over a Gaussian obstacle; 2D test case. The bathymetry is represented in grey.}
\label{fig:gausssoliton}
\end{figure}

In Figure \ref{fig:gausssoliton} the free surface elevation is represented, together with the bottom topography, at different times. The soliton is initially positioned in front of the obstacle (Figure \ref{fig:gausssolitona}). As soon as the soliton reaches the obstacle, its amplitude grows (Figure \ref{fig:gausssolitonb}). After the interaction with the obstacle the amplitude of the soliton gradually goes back to the initial value and a train of dispersive waves (Figures \ref{fig:gausssolitonc}, \ref{fig:gausssolitond}) can be observed behind the soliton.

\begin{figure}
\centering
\includegraphics[width=0.5\textwidth,trim=10 10 10 10,clip]{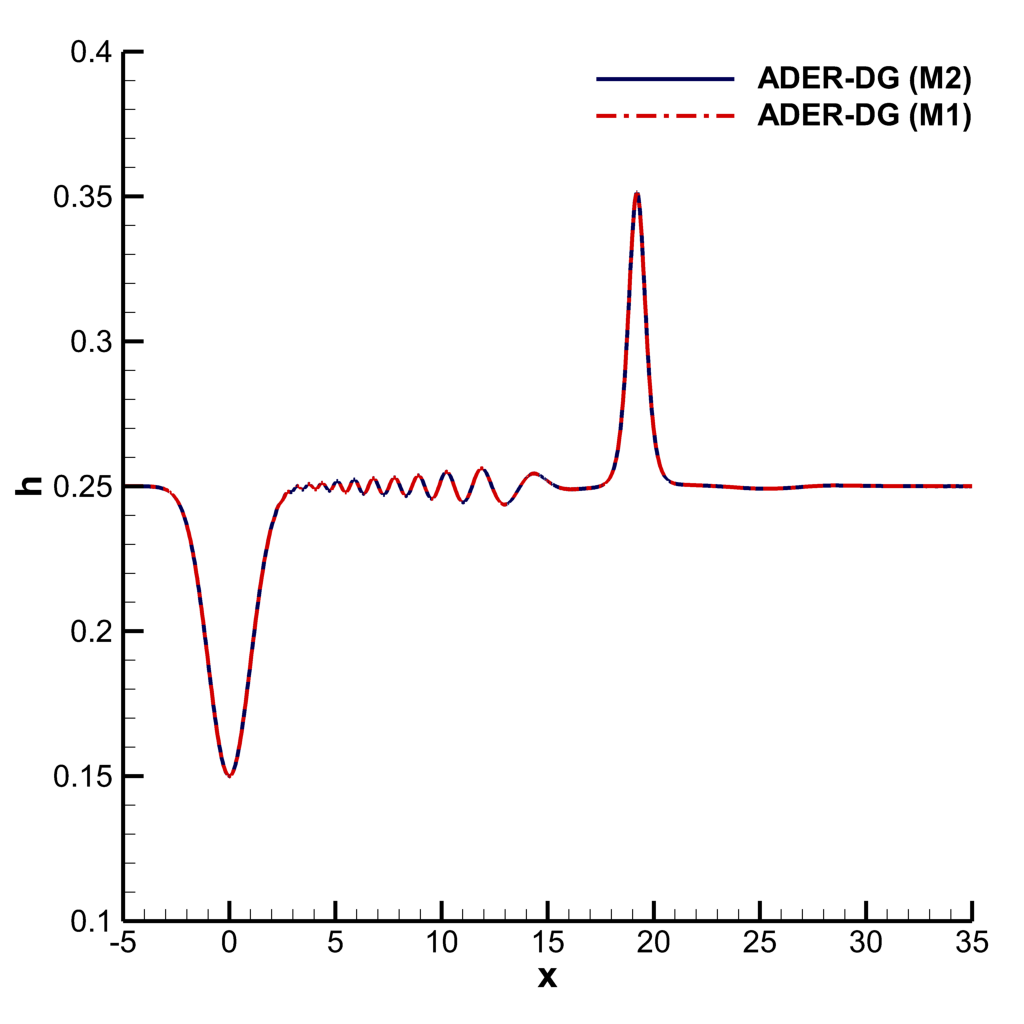}
\caption{One-dimensional cross section of the water depth, $h$, for $y=0$ at time $t_{\mathrm{end}}=12$, for the two-dimensional soliton over a Gaussian obstacle test case. Mesh $M1$: simulation with $N_x=N_y=200$ spatial elements. Mesh $M2$: simulation with $N_x=N_y=400$ spatial elements.}
\label{fig:gaussoliton_conv}
\end{figure}

In order to verify mesh convergence, a simulation with a refined spatial grid made of $400$ elements in both $x$ and $y$ directions has also been realized. In Figure \ref{fig:gaussoliton_conv}, a one-dimensional cross section at $y=0$ of the water depth $h$ is represented at the final time $t_{\mathrm{end}}=12$. As we can, see both mesh resolutions provide almost identical results, which indicates that the problem is well resolved.

\section{Conclusions}
\label{sec:conclusion}
A novel hyperbolic reformulation of the Serre-Green-Naghdi model for fully nonlinear and weakly dispersive water waves has been proposed without mild bottom assumption. The new hyperbolic model for dispersive water waves has been successfully tested in the context of high order accurate ADER-DG schemes and has been validated against quasi exact, numerical and experimental reference solutions. 

The model proposed in this paper proves to have good mathematical properties. First, it consists of a first order hyperbolic system,  
so that the severe time-step restrictions due to higher order derivatives of standard Boussinesq-type models are eliminated, thus allowing to effectively employ high order accurate explicit time integration techniques.
From the modelling perspective, the proposed system allows to recover the original SGN system in the  incompressible limit, when the artificial sound speed $c$ tends to infinity. Moreover, an additional energy conservation law, which reduces to the energy conservation law of the original SGN system in the incompressible limit, can be formulated. 

Concerning the numerical experiments, the convergence test for the propagation of a soliton over a flat bottom provides a validation of both the new model proposed in this paper and for the employed numerical approach. 
Since the new model is expected to be able to deal also with steep bottom bathymetry, different numerical tests with uneven topographies have been realised and compared with experimental data.    
In all the performed tests, the new model provides results that are in good agreement with the experimental data and other numerical reference solutions that can be found in the literature. Moreover, for very steep  topographies, like for the propagation of a soliton over a step, the new hyperbolic reformulation of the complete SGN model provides better results than the hyperbolic version of the SGN model with the mild bottom approximation introduced in \cite{escalante:2018}.

Future work will concern a rigorous derivation of the present first order hyperbolic model from variational principles, similar to the approach introduced by Gavrilyuk et al. in \cite{HGN,Dhaouadi2018} for first order hyperbolic approximations of nonlinear dispersive systems, but including  also a spatially variable topography without mild bottom assumption.

\section*{Acknowledgements}
C.B. and S.B.U. have been supported throughout the work on this paper by post-doctoral grants financed by INDAM in the framework of the \textit{Progetto
Premiale FOE 2014 - Strategic Initiatives for the Environment and Security (SIES)}. All the authors would like to thank prof. V. Vespri for getting them involved in this project. 

M.D. acknowledges the financial support received from the Italian Ministry of Education, University and Research (MIUR) in the frame of the Departments of Excellence Initiative 2018--2022 attributed to DICAM of the University of Trento (grant L. 232/2016) and in the frame of the PRIN 2017 project. MD has also received funding from the University of Trento via the  \textit{Strategic Initiative Modeling and Simulation}. 

All authors are member of the INdAM GNCS group. 

\bibliographystyle{plain}
\bibliography{nhsw,References}

\end{document}